\documentclass[leqno]{aoscmsize11}%
\usepackage{amsfonts}
\usepackage{amssymb}
\usepackage{amsthm}
\usepackage{latexsym}
\usepackage{graphicx}
\usepackage{amsmath}%
\newtheorem{theorem}{Theorem}

\newtheorem{lemma}{Lemma}

\newtheorem{remark}{Remark}
\newtheorem{assumption}{Assumption}
\def\theequation{\arabic{equation}}
\NONUMBIB
\def\baselinestretch{1.5}
\begin{document}

\SPECFNSYMBOL{}{1}{2}{}{}{}{}{}{}
%
\AOSMAKETITLE

\AOSAMS{Primary 62F12, 62M05; secondary 60H10, 60J60.}
\AOSKeywords{{Diffusions with jumps, efficiency, discrete sampling}\\}
\AOStitle{VOLATILITY ESTIMATORS FOR DISCRETELY SAMPLED LEVY PROCESSES}
\AOSauthor{Yacine A\"{\i}t-Sahalia\footnote
{Supported in part by NSF Grant SES-0350772.}
and Jean Jacod\footnote{Supported in part by the CNRS.}}
\AOSaffil{Princeton University and Universit\'e de Paris-6}
\AOSlrh{YACINE AIT-SAHALIA AND JEAN JACOD}
\AOSrrh{VOLATILITY ESTIMATORS FOR LEVY PROCESSES}
\renewcommand{\baselinestretch}{1.0}
\AOSAbstract
{This paper provides rate-efficient estimators of the volatility parameter
in the presence of L\'evy jumps.}
\maketitle
%

\BACKTONORMALFOOTNOTE{1}
%
%

\renewcommand{\baselinestretch}{1.5}%

\section{Introduction.\label{sec:intro}}

In this paper, we continue the study started in \cite{yacjacod04a}, about the
estimation of parameters when one observes a L\'{e}vy process $X$ at $n$
regularly spaced times $\Delta_{n},2\Delta_{n},\ldots,n\Delta_{n}$, with
$\Delta_{n}$ going to $0$ as $n\rightarrow\infty$. In our earlier paper, we
were concerned with the asymptotic behavior of the Fisher information, with
the objective of establishing a benchmark for what efficient estimators are
able to achieve in that context. Now, we wish to exhibit estimators which both
achieve that rate and can be explicitly computed.

We want to estimate a positive parameter $\sigma,$ which we call volatility,
in the model
\begin{equation}
X_{t}=\sigma W_{t}+Y_{t},\label{eq:X=sW+aY}%
\end{equation}
where $W$ is a standard Wiener process or, more generally, a symmetric stable
process of index $\beta$, and the process $Y$ is another L\'{e}vy process
without Wiener part and with jumps \textquotedblleft
dominated\textquotedblright\ in a sense we make precise below by those of $W$.
Allowing for jumps is of great interest in mathematical finance, in the
diverse contexts of option pricing, testing for the presence of jumps in asset
prices, interest rate modelling, risk management, optimal portfolio choice,
stochastic volatility modelling or for the purpose of better describing asset
returns data (see the references cited in \cite{yacjacod04a}).

Our aim is to construct estimators for $\sigma$ which behave under the model
(\ref{eq:X=sW+aY}) \textquotedblleft as well as\textquotedblright\ under the
model
\begin{equation}
X_{t}=\sigma W_{t},\label{eq:X=sW}%
\end{equation}
asymptotically as $\Delta_{n}\rightarrow0$ and $n\rightarrow\infty$. This is
in line with the results of \cite{yacjacod04a}, in which we proved that
property for the Fisher information. In other words, we want to be able to
estimate the volatility parameter $\sigma$ at the same rate when $Y,$ a jump
perturbation of $W,$ is present as when it is not. In some applications, $Y$
may represent frictions that are due to the mechanics of the trading process,
or in the case of compound Poisson jumps it may represent the infrequent
arrival of relevant information related to the asset. Given that both $W$ and
$Y$ contribute to the overall observed noise in $X,$ it is not a priori
obvious that it should be possible to estimate $\sigma$ equally well (at least
in the rate sense) with and without $Y.$ Beyond the robustness to
misspecification risk that such a result affords, it also for instance paves
the way for risk management or option hedging that is able to target the
\textquotedblleft$W$ risk\textquotedblright\ (continuous when $\beta=2$)
separately from the \textquotedblleft$Y$ risk\textquotedblright%
\ (discontinuous).

We distinguish between a \textit{parametric} case, where the law of $Y\ $is
known, and a \textit{semiparametric} case, where it is not. We show that, in
the parametric case, one can find estimators which are asymptotically
efficient in the Cramer--Rao sense, meaning that the asymptotic estimation
variance is equivalent as $n\rightarrow\infty$ to the inverse of the Fisher
information for the model (\ref{eq:X=sW}) without the perturbation $Y$. This
is possible when the law of $Y$ is completely known. In the semiparametric
case, where that law is unknown, obtaining asymptotically efficient estimators
requires $\Delta_{n}$ to go fast enough to $0$; but we can then exhibit
estimators that are efficient uniformly when the law of $Y$ stays in a set
sufficiently separated from the law of $W$. And in general we can exhibit a
large class of estimators which are consistent and achieve a specified rate
(although not the efficient rate).

A distinctive feature of the present paper is that we construct estimators
which are as simple as possible to implement. For example, in the parametric
situation where the law of $Y$ is known, one can in principle compute the MLE,
which is of course efficient. In practice, this is hardly feasible, as the
likelihood function derived from the convolution of the densities of $W$ and
$Y$ will in most situations not be available in closed form. So we provide a
number of other -- much simpler -- estimators which are not as good (in the
sense of not reaching the Cramer-Rao lower bound in general)\ but not too bad
either (in the sense of achieving the efficient rate of convergence).

The paper is organized as follows. In Section \ref{sec:setting}, we specify
our estimating setting. Section \ref{sec:ESTEQ} is devoted to estimating
equations: the estimators we propose all fall in that class and we state a
general result which covers them all. Sections \ref{sec:ParEst} and
\ref{sec:SemiParEst} are devoted to the parametric and semiparametric cases
respectively. Some examples are developed in Section \ref{sec:Example1},
\ref{sec:Example2}, \ref{sec:Example3} and \ref{sec:Example4}, where we
consider specific types of estimating equations such as the empirical
characteristic function , power variations and power variations with truncation.

\section{The setting.\label{sec:setting}}

With $X_{0}=0$, we observe $n$ i.i.d. increments from the L\'{e}vy process
(\ref{eq:X=sW+aY}),
\begin{equation}
\chi_{i}^{n}=X_{i\Delta_{n}}-X_{(i-1)\Delta_{n}}. \label{eq:Q1}%
\end{equation}
$W$ is a symmetric stable process of index $\beta\in(0,2]$, characterized by
\begin{equation}
\mathbb{E}(e^{iuW_{t}})=e^{-t|u|^{\beta}/2} \label{eq:Q2}%
\end{equation}
so that, when $\beta=2,$ $W$ is a standard Wiener process. The parameter to be
estimated is $\sigma$, and we will single out two situations concerning the
parameter space $\Theta$: either $\Theta=(0,\infty)$, or $\Theta$ is a compact
subset of $(0,\infty)$.

The law of $Y$ (as a process) is entirely specified by the law $G_{\Delta}$ of
the variable $Y_{\Delta}$ for any given $\Delta>0$. We write $G=G_{1}$, and we
recall that the characteristic function of $G_{\Delta}$ is given by the
L\'{e}vy-Khintchine formula
\begin{equation}
\mathbb{E}(e^{ivY_{\Delta}})=\exp\Delta\left(  ivb-{\frac{cv^{2}}{2}}+\int
F(dx)\left(  e^{ivx}-1-ivx1_{\{|x|\leq1\}}\right)  \right)  \label{eq:Ychar}%
\end{equation}
where $(b,c,F)$ is the \textquotedblleft characteristic
triple\textquotedblright\ of $G$ (or, of $Y$): $b\in\mathbb{R}$ is the drift
of $Y,$ and $c\geq0$ the local variance of the continuous part of $Y,$ and $F$
is the L\'{e}vy jump measure of $Y$, which satisfies $\int\left(  1\wedge
x^{2}\right)  F(dx)<\infty.$ We will denote by $\mathbf{P}_{\sigma,G}$ the law
of the process $X$.

We make $Y$ \textquotedblleft dominated\textquotedblright\ by $W$ in the
following sense: $G$ belongs to the class $\mathcal{G}_{\beta}$, defined as
follows. Let first $\Phi$ be the class of all increasing and bounded functions
$\phi:~(0,1]\rightarrow\mathbf{R}_{+}$ having $\lim_{x\downarrow0}\phi(x)=0$.
Then we set
\begin{align}
\mathcal{G}(\phi,\alpha)~  &  =\mbox{the set of all infinitely divisible
distributions with ~$c=0$~and,
for all $~x\in(0,1],$}\label{eq:Gclass}\\
&  \mbox{then }~\left\{
\begin{array}
[c]{ll}%
x^{\alpha}F([-x,x]^{c})\leq\phi(x)\qquad & \mbox{if }~\alpha<2\\
x^{2}F([-x,x]^{c})\leq\phi(x)\quad\mbox{and}~~\int_{\{|y|\leq x\}}%
|y|^{2}F(dy)\leq\phi(x)~\qquad & \mbox{if }~\alpha=2,
\end{array}
\right. \nonumber
\end{align}%
\begin{equation}
\mathcal{G}^{\prime}(\phi,\alpha)~=~\{G\in\mathcal{G}(\phi,\alpha
),~~G~~\mbox{is symmetrical about $0$}\}, \label{eq:G'class}%
\end{equation}%
\begin{equation}
\mathcal{G}_{\alpha}=\cup_{\phi\in\Phi}~\mathcal{G}(\phi,\alpha),\qquad
\mathcal{G}_{\alpha}^{\prime}=\cup_{\phi\in\Phi}~\mathcal{G}^{\prime}%
(\phi,\alpha), \label{eq:BGclass}%
\end{equation}
and we have
\begin{equation}
\left\{
\begin{array}
[c]{lll}%
\alpha\in(0,2]\quad & \Rightarrow\quad & \mathcal{G}_{\alpha}=\left\{
G~\text{is infinitely divisible},~c=0,~\lim_{x\downarrow0}x^{\alpha
}F([-x,x]^{c})=0\right\} \\
\alpha=2\quad & \Rightarrow\quad & \mathcal{G}_{2}=\left\{  G~\text{is
infinitely divisible},~c=0\right\}  .
\end{array}
\right.  \label{eq:BGclass2}%
\end{equation}
\smallskip

Now we recall some results from \cite{yacjacod04a}. The variable $W_{1}$
admits a $C^{\infty}$ density $h_{\beta}$, which is differentiable in the
state variable (the derivative is denoted by $h_{\beta}^{\prime}$). Then we
set
\begin{equation}
\breve{h}_{\beta}(w)=h_{\beta}(w)+wh_{\beta}^{\prime}(w),\quad\widetilde
{h}_{\beta}(w)={\frac{\breve{h}_{\beta}(w)^{2}}{h_{\beta}(w)}},\quad
\overline{h}_{\beta}(w)=\frac{wh_{\beta}^{\prime}(w)}{h_{\beta}(w)}%
,\quad\mathcal{I}(\beta)=\int\widetilde{h}_{\beta}(w)dw, \label{eq:Icallig}%
\end{equation}
so in fact $\mathcal{I}(\beta)$ is the Fisher information when we estimate
$\sigma$ on the basis of the single observation $\sigma W_{1}$ and for the
parameter value $\sigma=1$. The functions $\breve{h}_{\beta}$ and
$\widetilde{h}_{\beta}$ and $\overline{h}_{\beta}$ are also $C^{\infty}$, and
satisfy for some constant $c_{\beta}$:
\begin{equation}
\left\{
\begin{array}
[c]{l}%
\beta<2\quad\Rightarrow\quad h_{\beta}(w)+|\breve{h}_{\beta}(w)|+|\widetilde
{h}_{\beta}(x)|\leq\frac{c_{\beta}}{1+|w|^{1+\beta}},\quad|\overline{h}%
_{\beta}(w)|\leq c_{\beta},\\
\beta=2\quad\Rightarrow\quad\breve{h}_{\beta}(w)=(1-w^{2})h_{\beta}%
(w),\quad\widetilde{h}_{\beta}(x)=(1-w^{2})^{2}~h_{\beta}(w),\quad\overline
{h}_{\beta}(w)=-w^{2},
\end{array}
\right.  \label{eq:Maj}%
\end{equation}
and of course $h_{2}(w)=e^{-w^{2}/2}/\sqrt{2\pi}$, so in particular
$\mathcal{I}(\beta)=2$.

If we have a single observation $X_{\Delta}$ there is a (finite) Fisher
information for estimating $\sigma$, which we denote by $I_{\Delta}(\sigma
,G)$. With $n$ observed increments the corresponding Fisher information
becomes
\begin{equation}
I_{n,\Delta_{n}}(\sigma,G)=nI_{\Delta_{n}}(\sigma,G).\label{eq:n-Fisher}%
\end{equation}
The main result of \cite{yacjacod04a}, as far as the parameter $\sigma$ is
concerned, is summarized in the following:

\begin{theorem}
\label{theorem:sigma} a) If $G\in\mathcal{G}_{\beta}$ we have as
$\Delta\rightarrow0$:
\begin{equation}
I_{\Delta}(\sigma,G)\rightarrow\frac{1}{\sigma^{2}}\mathcal{I}(\beta).
\label{FI-IN}%
\end{equation}

b) For any $\phi\in\Phi$ we have as $\Delta\rightarrow0$:
\begin{equation}
\sup_{G\in\mathcal{G}(\phi,\beta)}\left|  I_{\Delta}(\sigma,G)-{\frac
{\mathcal{I}(\beta)}{\sigma^{2}}}\right|  \rightarrow0. \label{eq:Iss_uniform}%
\end{equation}

c) For each $n$ let $G^{n}$ be the standard symmetric stable law of index
$\alpha_{n}$, with $\alpha_{n}$ a sequence strictly increasing to $\beta$.
Then for any sequence $\Delta_{n}\rightarrow0$ such that $(\beta-\alpha
_{n})\log\Delta_{n}\rightarrow0$ (i.e. the rate at which $\Delta
_{n}\rightarrow0$ is slow enough), the sequence of numbers $I_{\Delta_{n}}
(\sigma,G^{n})$ converges to a limit which is strictly less than
$\mathcal{I}(\beta)/\sigma^{2}$.
\end{theorem}

Part (a) of the above theorem and (\ref{eq:n-Fisher}) hint towards the
existence of estimators $\widehat{\sigma}_{n}$ such that $\sqrt{n}%
~(\widehat{\sigma}_{n}-\sigma)$ converges to a centered Gaussian variable with
variance $\sigma^{2}/\mathcal{I}(\beta)$ under $\mathbf{P}_{\sigma,G}$, when
$G\in\mathcal{G}_{\beta}$ is known: this is the parametric situation, and we
will propose such estimators in Section \ref{sec:ParEst} below. In the
semiparametric situation where $G$ is unknown, (c) suggests that we cannot
achieve the same rate, unless, as given in (b), we know that $G$ is in the
class $\mathcal{G}(\phi,\alpha)$ for some $\alpha<\beta$ and some function
$\phi\in\Phi$.

As a matter of fact, we can do slightly better. If $\phi(x)=\zeta>0$ for all
$x$, we can still define $\mathcal{G}(\phi,\alpha)$ by (\ref{eq:Gclass}),
although $\phi$ no longer belongs to $\Phi$. We denote such a class by
$\overline{\mathcal{G}}(\zeta,\alpha)$, that is we introduce the notation (we
do not need to distinguish $\alpha<2$ and $\alpha=2$ here):
\begin{equation}
\label{eq:GBclass}%
\begin{array}
[c]{ll}%
\overline{\mathcal{G}}(\zeta,\alpha)~= &
\mbox{the set of all infinitely divisible
distributions with ~$c=0$~and,
for all $~x\in(0,1],$}\\
& \mbox{then }~x^{\alpha}F([-x,x]^{c})\leq\zeta,
\end{array}
\end{equation}
\begin{equation}
\overline{\mathcal{G}}^{\prime}(\zeta,\alpha)~=~ \{G\in\overline{\mathcal{G}%
}(\zeta,\alpha),~~G~~ \mbox{is symmetrical about $0$}\}, \label{eq:G'Bclass}%
\end{equation}
\begin{equation}
\overline{\mathcal{G}}_{\alpha}=\cup_{\zeta>0}~ \overline{\mathcal{G}}%
(\zeta,\alpha),\qquad\overline{\mathcal{G}}^{\prime}_{\alpha}=\cup_{\zeta>0}~
\overline{\mathcal{G}}^{\prime}(\zeta,\alpha). \label{eq:BGBclass}%
\end{equation}
The connection with the previous classes is as follows:
\begin{equation}
\mathcal{G}(\phi,\alpha)\subset\overline{\mathcal{G}}(\phi(1),\alpha),
\qquad\mathcal{G}_{\alpha}\subset\overline{\mathcal{G}}_{\alpha} \subset
\cap_{\alpha^{\prime}>\alpha}\mathcal{G}_{\alpha^{\prime}},\qquad
\mathcal{G}_{2}=\overline{\mathcal{G}}_{2}. \label{eq:BGBclass2}%
\end{equation}
For example, $\mathcal{G}_{0}$ is the class of all $G$'s for which $Y$ is a
pure drift ($Y_{t}=bt$), whereas $\overline{\mathcal{G}}_{0}$ is the class of
all $G$'s for which $Y$ is a compound Poisson process plus a drift. Also, any
stable process $Y$ with index $\alpha<2$ belongs to $\overline{\mathcal{G}%
}_{\alpha}$, but not to $\mathcal{G}_{\alpha}$.

\section{About estimating equations.\label{sec:ESTEQ}}

The practical estimators we will propose for $\sigma$ are all obtained by
setting an estimating equation (also known as a generalized moment condition)
to zero. We prove here a general result about the asymptotic properties of
such estimators, which will be used several times below. Similar general
results for estimating equations are of course known (see various forms in
\cite{godambe60}, \cite{hansengmm82} and \cite{heyde97}), but we adapt them
here to our setting with assumptions (by no means minimal)\ that are
sufficient in our context.

Recall that we want to estimate a parameter $\sigma>0$. At stage $n$ we
observe $p_{n}$ i.i.d. random variables $\chi_{i}^{n}$ and introduce two
auxiliary variables $S_{n}>0$ and $Q_{n}\in\mathbf{R}$. Under the associated
probability measure $\mathbf{P}_{n,\sigma}$ we suppose that the families
$(S_{n},Q_{n})$ and $(\chi_{i}^{n}:1\leq i\leq p_{n})$ are independent, and of
course $p_{n}\rightarrow\infty$. Let us introduce the following conditions:

\begin{assumption}
[A1]If $\sigma_{n}\rightarrow\sigma>0$ then $S_{n}\rightarrow\sigma$ in
$\mathbf{P}_{n,\sigma_{n}}$--probability.
\end{assumption}

\begin{assumption}
[A2]If $\sigma_{n}\rightarrow\sigma>0$ then the sequence $(Q_{n}\mid
\mathbf{P}_{n,\sigma_{n}})$ is tight.
\end{assumption}

Next we consider two families $(f_{n,s,q})_{s>0}$ and $(H_{n,s})_{s>0,q\in
\mathbf{R}}$ of functions on $\mathbf{R}$ and $(0,\infty)$ respectively, to be
specified later but with adequate integrability and smoothness properties, and
we associate the estimating function\
\begin{equation}
U_{n,s,q}(u)=\frac{1}{p_{n}}\sum_{i=1}^{p_{n}}\left(  f_{n,s,q}(\chi_{i}%
^{n})-H_{n,s}(u)\right)  .\label{Q12}%
\end{equation}
In this exactly-identified context, we set
\begin{equation}
\widehat{\sigma}_{n}(s,q)=\left\{
\begin{array}
[c]{ll}%
\mbox{the $u>0$ with $U_{n,s,q}(u)=0$ which is closest to $s$}\quad &
\mbox{if it exists}\\
1 & \mbox{otherwise}
\end{array}
\right.  \label{Q13}%
\end{equation}
(if $U_{n,s,q}=0$ has two closest solutions at equal distance of $s$, we
select the smallest one). We also set
\begin{equation}
F_{n,s,q}(\sigma)=\mathbf{E}_{n,\sigma}(f_{n,s,q}(\chi_{i}^{n})),\qquad
F_{n,s,q}^{(2)}(\sigma)=\mathbf{E}_{n,\sigma}(f_{n,s,q}(\chi_{i}^{n}%
)^{2}).\label{Q11}%
\end{equation}
Note in particular that we are not assuming that the estimating equation is
correctly centered: correct centering would requiring using $F_{n,s,q}$
instead of $H_{n,s}.$ $H_{n,s}$ may be equal to $F_{n,s,q},$ but can also be
just an approximation to it (in which case we will talk about
\textquotedblleft approximate centering\textquotedblright) that may for
instance be valid as $n\rightarrow\infty$. Incorrect centering leads to
estimators that are asymptotically biased, although that effect can be
mitigated as $n\rightarrow\infty$ if $H_{n,s}$ approximates $F_{n,s,q}$ (see
Assumption (B5) below).

Let us now list a series of assumptions on the previous functions:

\begin{assumption}
[B1]We have $\sup_{n\geq1,~s>0,~q\in\mathbf{R}}\Vert f_{n,s,q}\Vert^{4}%
/p_{n}<\infty$, where $\Vert f\Vert$ is the sup--norm.
\end{assumption}

\begin{assumption}
[B2]$H_{n,s}$ is continuously differentiable.
\end{assumption}

\begin{assumption}
[B3]For all $s>0$ there is a differentiable function $\overline{F}_{s}$ on
$(0,\infty)$, such that whenever $s_{n}\rightarrow s$ then $H_{n,,s_{n}}$ and
$H_{n,s_{n}}^{\prime}$ converge locally uniformly to $\overline{F}_{s}$ and
$\overline{F}_{s}^{\prime}$ respectively.
\end{assumption}

\begin{assumption}
[B4]$\overline{F}_{s}^{\prime}(s)\neq0$ for all $s>0$.
\end{assumption}

\begin{assumption}
[B5]$F_{n,s_{n},q_{n}}^{(2)}(u_{n})$ converges to a limit $F^{(2)}(u)$ for any
two sequences $u_{n}$ and $s_{n}$ converging to the same limit $u>0$ and any
bounded sequence $q_{n}$.
\end{assumption}

\begin{assumption}
[B6] There is a sequence $w_{n}\rightarrow+\infty$ such
that $\sup_{n}~w_{n}|F_{n,s_{n},q_{n}}(u_{n})-H_{n,s_{n}}(u_{n}))|<\infty$ for
any two sequences $u_{n}$ and $s_{n}$ converging to the same limit $u>0$ and
any bounded sequence $q_{n}.$
\end{assumption}

Then we have the following:

\begin{theorem}
\label{TTT1} Assume (A1), (A2) and (B1)--(B6).

a) The sequence $((w_{n}\bigwedge\sqrt{p_{n}})(\widehat{\sigma}_{n}%
(S_{n},Q_{n})-\sigma_{n}))$ is tight under $\mathbf{P}_{n,\sigma}$, uniformly
in $n$ and in $\sigma$ in any compact subset of $(0,\infty)$.

b) If $w_{n}/\sqrt{p_{n}}\to\infty$, then the sequence $(\sqrt{p_{n}}~
(\widehat{\sigma}_{n}(S_{n},Q_{n})-\sigma_{n}))$ converges in law under
$\mathbf{P}_{n,\sigma}$, uniformly in $\sigma$ in any compact subset of
$(0,\infty)$, towards the centered normal distribution with variance $\Xi
^{2}(\sigma):=\left(  F^{(2)}(\sigma)-\overline{F}_{\sigma}(\sigma
)^{2}\right)  /\overline{F}^{\prime}_{\sigma}(\sigma)^{2}$.
\end{theorem}

We devote the remainder of this section to proving this theorem. First, we
state a lemma which gathers some classical limit theorems on i.i.d. triangular
arrays. For each $n$ let $(\zeta_{i}^{n}:i=1,\ldots,q_{n})$ be real--valued
and i.i.d. random variables, possibly defined on different probability spaces
$(\Omega_{n},\mathcal{F}_{n},\mathbf{P}_{n})$ when $n$ varies. Then:

\begin{lemma}
\label{LE3} Assume that $\zeta^{n}_{i}$ is square--integrable, and set
$\gamma_{n}=\mathbf{E}_{n}(\zeta^{n}_{i})$ and $\Gamma_{n}=\mathbf{E}%
_{n}((\zeta^{n}_{i})^{2})-\gamma_{n}^{2}$. If $p_{n}\to\infty$ and $\Gamma
_{n}/p_{n}\to0$, we have
\begin{equation}
\label{LLN}\frac1{p_{n}}\sum_{i=1}^{p_{n}}\zeta^{n}_{i}-\gamma_{n}~
\overset{\mathbf{L}^{2}(\mathbf{P}_{n})}{\longrightarrow}~0.
\end{equation}
Furthermore if $\Gamma_{n}\to\Gamma$ for some limit $\Gamma>0$ and if
$\mathbf{E}(|\zeta^{n}_{i}|^{4})/p_{n}\to0$, we have
\begin{equation}
\label{CLT}\sqrt{p_{n}}~\left(  \frac1{p_{n}}\sum_{i=1}^{p_{n}}\zeta^{n}%
_{i}-\gamma_{n}\right)  ~ \overset{\mathcal{L}(\mathbf{P}_{n})}%
{\longrightarrow}~ \mathcal{N}(0,\Gamma).
\end{equation}

\end{lemma}

In the next three lemmas we suppose that $\sigma_{n}\rightarrow\sigma>0$, and
we write $\mathbf{P}_{n}=\mathbf{P}_{n,\sigma_{n}}$.

\begin{lemma}
\label{LE5} Let $s_{n}\to\sigma$ and let $q_{n}$ be a bounded sequence.

a) The sequence $\left(  (w_{n}\bigwedge\sqrt{p_{n}})~U_{n,s_{n},q_{n}}%
(\sigma_{n}) \mid\mathbf{P}_{n}\right)  $ is tight.

b) If $w_{n}/\sqrt{p_{n}}\to\infty$ then
\[
\sqrt{p_{n}}~U_{n,s_{n},q_{n}}(\sigma_{n})~ \overset{\mathcal{L}%
(\mathbf{P}_{n})}{\longrightarrow}~\mathcal{N}(0,F^{(2)}(\sigma)-\overline
{F}_{\sigma}(\sigma)^{2}).
\]

\end{lemma}

\begin{proof}
We have $U_{n,s_{n},q_{n}}(\sigma_{n}) =\frac1{p_{n}}\sum_{i=1}^{p_{n}}%
\zeta^{n}_{i}$, where for each $n$ the $\zeta^{n}_{i}$'s are i.i.d. with mean
and variance given by
\[
\gamma_{n}=F_{n,s_{n},q_{n}}(\sigma_{n})-H_{n,s_{n}}(\sigma_{n}),\qquad
\Gamma_{n}=F^{(2)}_{n,s_{n},q_{n}}(\sigma_{n})-F_{n,s_{n},q_{n}}(\sigma
_{n})^{2},
\]
and further $|\zeta^{n}_{i}|\leq\alpha_{n}$ for numbers $\alpha_{n}$
satisfying $\alpha_{n}^{4}/p_{n}\to0$ by (B1). Now (B6) yields that
$\gamma_{n}\to0$, hence (B3) yields $F_{n,s_{n},q_{n}}(\sigma_{n})\to
\overline{F}_{\sigma}(\sigma)$. On the other hand, (B5) implies $F^{(2)}%
_{n,s_{n},q_{n}}(\sigma_{n})\to F^{(2)}(\sigma)$.

Therefore it follows from (\ref{CLT}) that
\begin{equation}
\label{eq:CLT1}\sqrt{p_{n}}~\left(  U_{n,s_{n},q_{n}}(\sigma_{n})-\gamma
_{a}\right)  ~ \overset{\mathcal{L}(\mathbf{P}_{n})}{\longrightarrow
}~\mathcal{N}(0,F^{(2)}(\sigma)-\overline{F}_{\sigma}(\sigma)^{2}),
\end{equation}
and since $\sup_{n}~w_{n}|\gamma_{n}|<\infty$ by (B6), we readily get
the two results.
\end{proof}

\begin{lemma}
\label{LE6} a) The sequence $((w_{n}\bigwedge\sqrt{p_{n}})~U_{n,S_{n},Q_{n}%
}(\sigma_{n})\mid\mathbf{P}_{n})$ is tight.

b) If $w_{n}/\sqrt{p_{n}}\to\infty$, the sequence $(\sqrt{p_{n}}~
U_{n,S_{n},Q_{n}}(\sigma_{n})\mid\mathbf{P}_{n})$ converges in law towards the
centered normal distribution with variance $F^{(2)}(\sigma)-\overline
{F}_{\sigma}(\sigma)^{2}$.
\end{lemma}

\begin{proof}
a) Let $V(n,s,q)=(w_{n}\bigwedge\sqrt{p_{n}})U_{n,,s,q}(\sigma_{n})$. The
previous lemma implies that as soon as the deterministic sequence $s_{n}$
converges to $\sigma$, we have for all $B>0$:
\begin{equation}
\label{eq:Ti}\lim_{A\to\infty}~\sup_{n\geq1}~u_{A,B}(n,s_{n}) =0,\qquad
\mbox{where }~ u_{A,B}(n,s)=\sup_{|q|\leq B}~\mathbf{P}_{n}(|V(n,s,q)|>A).
\end{equation}
If the sequence $(V(n,S_{n},Q_{n}) \mid\mathbf{P}_{n})$ is not tight, there
exists an infinite sequence $n_{k}$ such that $\mathbf{P}_{n_{k}}%
(|V(n_{k},S_{n_{k}},Q_{n_{k}})|>A)\geq1/A$ for some $A>0$ and, up to taking a
further subsequence still denoted by $n_{k}$ we can assume by (A1) that
$S_{n_{k}}\to\sigma$ pointwise. Since $(S_{n},Q_{n})$ is independent of the
family $(V(n,s,q);s>0,q\in\mathbf{R})$, we get
\[
\mathbf{P}_{n_{k}}(|V(n_{k},S_{n_{k}},Q_{n_{k}})|>A)\leq\mathbf{P}_{n_{k}%
}(|Q_{n_{k}}|>B)+\mathbf{E}_{n_{k}}(u_{A,B}(n_{k},S_{n_{k}}))).
\]
Then (\ref{eq:Ti}) and Lebesgue's Theorem imply that
\[
\limsup_{k}~\mathbf{P}_{n_{k}}(|V(n_{k},S_{n_{k}},Q_{n_{k}})|>A)\leq\sup
_{n}~\mathbf{P}_{n}(|Q_{n}|>B)
\]
for all $B>0$ and, in view of (A2), we deduce that $\limsup_{k}~\mathbf{P}%
_{n_{k}}(|V(n_{k},S_{n_{k}},Q_{n_{k}})|>A)=0$: this contradicts the definition
of the sequence $n_{k}$, and we have the result.

b) Let us denote by $V$ a variable with law $\nu=\mathcal{N}(0,F^{(2)}%
(\sigma)-\overline{F}_{\sigma}(\sigma))$. Let $\nu_{n,s,q}$ be the law of
$V(n,s,q):=\sqrt{p_{n}}~U_{n,s,q}(\sigma_{n})$. The claim amounts to proving
that, for all bounded continuous functions $g$, we have
\begin{equation}
\mathbf{E}_{n}\left(  g(V(n,S_{n},Q_{n}))\right)  ~\rightarrow~\mathbf{E}%
(g(V)).\label{eq:CLT3}%
\end{equation}
For this, it is enough to prove that from any subsequence one can extract a
further subsequence along which (\ref{eq:CLT3}) holds. So, in view of (A1) and
(A2) it is no restriction to assume that in fact $(S_{n},Q_{n})$ converges in
law to $(\sigma,Q)$ for some variable $Q$.

In fact, due to the independence of $(S_{n},Q_{n})$ and $(W^{\prime
}(n,s,q):s>0,q\in\mathbf{R})$, we can replace the pair $(S_{n},Q_{n})$ in the
left side of (\ref{eq:CLT3}) by any other pair $(S_{n}^{\prime},Q_{n}^{\prime
})$ having the same law than $(S_{n},Q_{n})$ and still independent of
$(W^{\prime}(n,s,q):s>0,q\in\mathbf{R})$. Therefore, using the Skorokhod
representation theorem, we can indeed assume that $(S_{n},Q_{n})$ converges
pointwise to $(\sigma,Q)$. Then
\[
\mathbf{E}_{n}\left(  g(V(n,S_{n},Q_{n}))\right)  =\mathbf{E}_{n}\left(
\int\nu_{n,S_{n},Q_{n}}(dx)g(x))\right)  .
\]
Since $S_{n}\rightarrow\sigma$ and $Q_{n}\rightarrow Q$, one deduces from
Lemma \ref{LE5}--(b) that the sequence $\int\nu_{n,S_{n},Q_{n}}(dx)g(x)$
converges pointwise to $\int\nu(dx)g(x)=\mathbf{E}(g(V))$, and it is bounded
by $\Vert g\Vert$, so Lebesgue's Theorem yields (\ref{eq:CLT3}).
\end{proof}

\begin{lemma}
\label{LE4} The sequence $\widehat{\sigma}_{n}$ converges in $\mathbf{P}_{n}%
$--probability to $\sigma$.
\end{lemma}

\begin{proof}
Exactly as in the previous proof, without loss of generality we can assume
that the pair $(S_{n},Q_{n})$ converges pointwise to $(\sigma,Q)$ with $Q$ a
suitable random variable.

Lemma \ref{LE6} implies that $U_{n,S_{n},Q_{n}}(\sigma_{n}) \to0$ in
probability (recall that both $w_{n}$ and $p_{n}$ go to infinity). Observe
that
\[
U_{n,S_{n},Q_{n}}(u)-U_{n,S_{n},Q_{n}}(\sigma_{n})= H_{n,S_{n}}(\sigma
_{n})-H_{n,S_{n}}(u),
\]
which by (B3) converges (pointwise) locally uniformly in $u$ towards
$H(u):=\overline{F}_{\sigma}(\sigma)-\overline{F}_{\sigma}(u)$. Hence
$U_{n,S_{n},Q_{n}}(u)$ also converges locally uniformly in $u$ towards $H(u)$,
in $\mathbf{P}_{n}$--probability. But by (B4) the function $H$ is null at
$\sigma$ and is either strictly decreasing or strictly increasing in a
neighborhood of $\sigma$: then the definition (\ref{Q13}) of $\widehat{\sigma
}_{n}(S_{n},Q_{n})$ immediately gives the result.
\end{proof}

Finally, we prove Theorem \ref{TTT1}:

\begin{proof}
[Proof of Theorem \ref{TTT1}]As usual, to get the local uniformity in $\sigma$
for the tightness in (a) or the convergence in (b), it is enough to obtain the
tightness (resp. convergence) under $P_{n}=P_{n,\sigma_{n}}$ for any sequence
$\sigma_{n}\rightarrow\sigma>0$. Let us write for simplicity $\widehat{\sigma
}_{n}=\widehat{\sigma}_{n}(S_{n},Q_{n})$ and $U_{n}=U_{n,S_{n},Q_{n}}$.

By (B2), $U_{n}$ is continuously differentiable. We deduce from Lemma
\ref{LE4} the existence of sets $A_{n}$ with $P_{n}(A_{n})\rightarrow1$, such
that on $A_{n}$ we have $U_{n}^{\prime}(\widehat{\sigma}_{n})=0$, and thus
Taylor's formula yields a random variable $T_{n}$ taking its values between
$\sigma_{n}$ and $\widehat{\sigma}_{n}$, and such that
\begin{equation}
U_{n}(\sigma_{n})=-(\widehat{\sigma}_{n}-\sigma_{n})U_{n}^{\prime}%
(T_{n})\qquad\mbox{on the set }~A_{n}. \label{eq:Q20}%
\end{equation}

Observe that $U_{n}^{\prime}(T_{n})=-H_{n,S_{n}}^{\prime}(T_{n})$,. Since both
$S_{n}$ and $T_{n}$ converge in probability to $\sigma$, (B3) implies that
$U_{n}^{\prime}(T_{n})\rightarrow-\overline{F}_{\sigma}^{\prime}(\sigma)$ in
probability. Since $\overline{F}_{\sigma}^{\prime}(\sigma)\neq0$ by (B4), all
the results of our theorem are now easily deduced from (\ref{eq:Q20}) and
Lemma \ref{LE6}.
\end{proof}

With this general result in hand, we now turn to our specific situation:
estimating $\sigma$ in the presence of the L\'{e}vy process $Y,$ first when
the law of $Y$ is known and second when it is not.

\section{Estimation of $\sigma$ in the parametric case.\label{sec:ParEst}}

In this section, we study the estimation of $\sigma$ when the law of $Y,$
i.e., the measure $G\in G_{\beta}$, is known. We will construct a class of
estimating equations for $\sigma$, with $\chi_{i}^{n}$ given by (\ref{eq:Q1}).

\subsection{Construction of the estimators.}

In the sequel the number $\beta\in(0,2]$ is fixed and does not usually appear
explicitly in our notation. A constant which depends only on $\beta$ and on
another parameter $\gamma$ is denoted by $C_{\gamma}$, and it may change from
line to line. If $G\in G_{\alpha}$ with $\alpha\leq\beta$, and with the
associated process $Y$, we set
\begin{equation}
b^{\prime}(G,\alpha)=\left\{
\begin{array}
[c]{ll}%
b-\int_{\{|x|\leq1\}}xF(dx)\text{ \ \ \ } & \text{if \ }\alpha<1\\[2.5mm]%
b & \text{if \ }\alpha\geq1,
\end{array}
\right.  \qquad Z_{\Delta}(\alpha):=\Delta^{-1/\beta}\left(  Y_{\Delta
}-b^{\prime}(G,\alpha)\Delta\right)  \label{eq:b'}%
\end{equation}
and we let $G_{\Delta,\alpha}^{\prime}$ denote the law of $Z_{\Delta}(\alpha
)$. Then we define the \textquotedblleft modified increments\textquotedblright%
\ (recall (\ref{eq:Q1})):
\begin{equation}
\chi_{i}^{\prime n}(G)=\Delta_{n}^{-1/\beta}(\chi_{i}^{n}-b^{\prime}%
(G,\beta)\Delta_{n}). \label{FE2}%
\end{equation}

Next, for any $\alpha\in(0,2]$ and any $\phi\in\Phi$ we set for $x\in(0,1)$:
\begin{equation}
\phi_{\alpha}(x)=\left\{
\begin{array}
[c]{ll}%
\frac{\phi(x)}{1-\alpha} & \mbox{if }~\alpha<1\\[2.5mm]%
\phi(x)+\frac{\phi(x)}{\sqrt{\log(1/x)}}+\phi\left(  1\wedge e^{-\sqrt
{\log(1/x)}}\right)  \qquad & \mbox{if }~\alpha=1\\[2.5mm]%
\phi(x)+\frac{\phi(\sqrt{x})}{\alpha-1}+\frac{\phi(1)}{\alpha-1}%
~x^{\frac{\alpha-1}{2}} & \mbox{if }~\alpha>1.
\end{array}
\right.  \label{eq:phialp}%
\end{equation}
This defines an increasing function $\phi_{\alpha}:~(0,1]\rightarrow R_{+}$
having $\phi\leq\phi_{\alpha}$ and $\phi_{\alpha}(x)\rightarrow0$ as
$x\rightarrow0$.

Next, if $G\in G_{\alpha}$ for some $\alpha\leq\beta$, and $u>0$ and $v\geq0$
and $z\in R$ and if $k$ is a bounded function, we set
\begin{equation}
\Psi_{G,\Delta,\alpha,k}(u,v,z)=\int h_{\beta}(x)dx\int G_{\Delta,\alpha
}^{\prime}(dw)~k(ux+vw+z). \label{FF1}%
\end{equation}

Finally, we introduce the \textquotedblleft tail function\textquotedblright\
\begin{equation}
\psi(u)=\mathbf{P}(|W_{1}|>1/u)=2\int_{1/u}^{\infty}h_{\beta}(x)dx \label{FE3}%
\end{equation}
for $u>0$ (this depends on $\beta$): it is $C^{\infty}$, strictly increasing
from $0$ to $1$, with non-vanishing first derivative. So its reciprocal
function $\psi^{-1}$, from $(0,1)$ into $(0,\infty)$, is also $C^{\infty}$ and
strictly increasing.

Recall that we work here under the assumption that $G\in G_{\beta}$ is known,
and so in particular we know $b^{\prime}(G,\beta)$; we also have $G\in
G(\phi,\beta)$ for some $\phi\in\Phi$. We need first a preliminary estimator,
which is constructed as follows. We choose an arbitrary sequence $m_{n}$ of
integers satisfying
\begin{equation}
m_{n}\uparrow\infty,\qquad\frac{m_{n}}{n}\rightarrow0 \label{FE1}%
\end{equation}
and, recalling (\ref{FE2}) and (\ref{FE3}), we set
\begin{equation}
V_{n}(G)=\frac{1}{m_{n}}~\sum_{i=1}^{m_{n}}1_{\{|\chi_{i}^{\prime n}%
(G)|>1\}},\qquad S_{n}(G)=\left\{
\begin{array}
[c]{ll}%
\psi^{-1}(V_{n}(G))\quad & \mbox{if }~0<V_{n}(G)<1\\[2mm]%
1 & \mbox{otherwise}.
\end{array}
\right.  \label{FE4}%
\end{equation}

To form an estimating equation for the construction of the final estimator of
$\sigma$, we choose a function $k$ satisfying
\begin{equation}
\sup_{x}~\frac{|k(x)|}{1+|x|^{\gamma}}<\infty,\qquad I(k):=\int\breve
{h}_{\beta}(x)k(x)dx\neq0, \label{eq:k}%
\end{equation}
where the number $\gamma$ satisfies
\begin{equation}
\gamma\geq0,\qquad\beta\leq2\quad\Rightarrow\quad\gamma<\frac{\beta}{2}.
\label{eq:ga}%
\end{equation}

Then we set
\begin{equation}
k_{n}(x)=\left\{
\begin{array}
[c]{ll}%
k(x) & \mbox{if $k$ is bounded}\\[2.5mm]%
k(x)~1_{\{|k(x)|\leq\nu_{n}\}}\qquad & \mbox{otherwise},
\end{array}
\right.  \label{eq:knn}%
\end{equation}
where $\nu_{n}$ be an increasing sequence of numbers satisfying
\begin{equation}
\nu_{n}\rightarrow\infty,\qquad\nu_{n}^{2}~\phi_{\beta}(\Delta_{n}^{1/\beta
})~\rightarrow~0,\qquad\frac{\nu_{n}^{4}}{n}~\rightarrow~0,\label{eq:kn}%
\end{equation}
and where $\phi_{\beta}$ is associated with $\phi$ (a function such that $G\in
G(\phi,\beta)$) by (\ref{eq:phialp}). Then, with the notation $p_{n}=n-m_{n}$,
and since each $k_{n}$ is bounded, we can define the following estimation
functions (for $u>0$):
\begin{equation}
U_{n,G,\phi,k}(u)=\frac{1}{p_{n}}\sum_{i=m_{n}+1}^{n}k_{n}\left(  \frac
{\chi_{i}^{\prime n}(G)}{S_{n}(G)}\right)  -\Psi_{G,\Delta_{n},\beta,k_{n}%
}\left(  \frac{u}{S_{n}(G)},\frac{1}{S_{n}(G)},0\right)  .\label{FE6}%
\end{equation}

Finally the estimators for $\sigma$ are:
\begin{equation}
\widehat{\sigma}_{n}(G,\phi,k)=\left\{
\begin{array}
[c]{ll}%
\mbox{the $u>0$~~with ~$U_{n,G,\phi,k}(u)=0$ ~which is closest to
$S_n(G)$}\quad & \mbox{if it exists}\\[2.5mm]%
1 & \mbox{otherwise}.
\end{array}
\right.  \label{FE7}%
\end{equation}
As the notation suggests, this estimator depend on $G$ and on $k$ in an
obvious way, and it depends on $\phi$ through the choice for $k_{n}$ made in
(\ref{eq:kn}). It also depends on $\beta$, but we leave this dependency
implicit to avoid cluttering the notation.

\subsection{Asymptotic distribution in the parametric case.}

With the function $k$ as in (\ref{eq:k}), the following defines two finite
numbers:
\begin{equation}
J(k)=\mathbf{E}(k(W_{1})^{2})-(\mathbf{E}(k(W_{1})))^{2},\qquad\Sigma
^{2}(k)=\frac{J(k)}{I(k)^{2}}. \label{eq:FF}%
\end{equation}

\begin{theorem}
\label{theorem:Gknown} Let $\phi\in\Phi$, and let $k$ be a function satisfying
(\ref{eq:k}) for some $\gamma$ having (\ref{eq:ga}). Suppose also that
$\Delta_{n}\rightarrow0$.

a) The sequence $\sqrt{n}~(\widehat{\sigma}_{n}(G,\phi,k)-\sigma)$ converges
in law to $N(0,\sigma^{2}\Sigma^{2}(k))$, under $P_{\sigma,G}$, uniformly in
$G\in G(\phi,\beta)$ and in $\sigma\in\lbrack\varepsilon,1/\varepsilon]$ for
any $\varepsilon>0$.

b) We have $\Sigma^{2}(k)\geq1/I(\beta)$, and this inequality is an equality
if we choose $k=\overline{h}_{\beta}$.
\end{theorem}

Now we give a number of comments and examples.

\begin{remark}
In light of (\ref{eq:FF}), it is of course possible / advisable to select the
function $k$ to minimize $\Sigma^{2}(k).$ The choice $k=\overline{h}_{\beta}$
is indeed possible: by (\ref{eq:Maj}) the function $k=\overline{h}_{\beta}$
satisfies (\ref{eq:k}) with $\gamma=0$ (resp. $\gamma=2$) if $\beta<2$ (resp.
$\beta=2$). Such a choice gives asymptotically efficient estimators, in the
strong\ sense that they behave asymptotically like the efficient estimators
for the model $X_{t}=\sigma W_{t}$ (with no perturbing term $Y$).
\end{remark}

\begin{remark}
To put these estimators in use we would need to numerically compute the
function $\Psi_{G,\Delta,\beta,k}(u,v,0)$, for a single value of $v$ (either
$1$ or $1/S_{n}(G)$), and all\ values of $u$ (in principle). Except in special
situations (see for instance Section \ref{sec:Example1}), there is no closed
form for this function, and we have to resort to numerical integration or to
Monte--Carlo techniques. For this it is of course helpful to have a closed
form for $k$ (or rather for the truncated $k_{n}$). In general, this is not
the case for the function $k=\overline{h}_{\beta}$ (the optimal choice),
unless $\beta=2$.
\end{remark}

\begin{remark}
As an example of function $k$, we can take $k(x)=|x|^{r}$, for some $r>0$ when
$\beta=2$ and $r\in(0,\beta/2)$ otherwise (when $\beta=2$ and $r=2$ this is
the optimal choice since $\overline{h}_{2}(x)=-x^{2}$): the function
$\Psi_{G,\Delta_{n},\beta,k_{n}}$ is still not explicit, but it is easily
approximated by Monte--Carlo techniques, at last when $Y_{t}$ can be
simulated, or it may be available in closed form for some common distributions
of $Y$. We will do that in some detail in Section \ref{sec:Example2}. In any
event, the limiting variance is easy to compute from (\ref{eq:FF}).
\end{remark}

\begin{remark}
Another possibility is to use the empirical characteristic function\ of the
sampled increments, which leads to an closed form expression for
$\Psi_{G,\Delta_{n},\beta,k_{n}}$. This will be done in Section
\ref{sec:Example1}.
\end{remark}

\subsection{Some preliminaries.}

Here we gather some results from \cite{yacjacod04a}, and also about the
functions of (\ref{FF1}), which will be used to obtain the previous theorem
and for further results as well. First we recall Lemma 2 of \cite{yacjacod04a}%
: for any $\phi\in\Phi$, and with the notation (\ref{eq:phialp}), we have for
$\Delta\leq1$ and $\alpha\leq\beta$ and $K\geq0$ and some constant
$C=C_{\alpha}$ depending on $\alpha$ only,
\begin{equation}
G\in\mathcal{G}(\phi,\alpha),~~|g(x)|\leq K(1\wedge|x|)\quad\Longrightarrow
\quad\mathbf{E}(|g(Z_{\Delta}(\alpha)|)\leq CK\Delta^{\frac{2(\beta-\alpha
)}{\beta(2+\alpha)}}\phi_{\alpha}(\Delta^{\frac{2+\beta}{\beta(2+\alpha)}}).
\label{eq:Z'}%
\end{equation}
In fact the proof of this result also works when $\phi(x)=\zeta$ for all $x$
(with $\phi_{\alpha}$ substituted with a constant), thus giving
\begin{equation}
G\in\overline{\mathcal{G}}(\zeta,\alpha),~~|g(x)|\leq K(1\wedge|x|)\quad
\Longrightarrow\quad\mathbf{E}(|g(Z_{\Delta}(\alpha)|)\leq CK\zeta
\Delta^{\frac{2(\beta-\alpha)}{\beta(2+\alpha)}}. \label{eq:Z'modif}%
\end{equation}

This is not enough for our purposes, at least in the semiparametric situation,
and we will need also the next lemma about symmetrical measures:

\begin{lemma}
\label{lem:1} If $\Delta\leq1$ and $\alpha\leq\beta$ and $K\geq0$, we have for
some constant $C$ depending on $\alpha$ only:
\begin{equation}
G\in\overline{\mathcal{G}}^{\prime}(\zeta,\alpha),~~|g(x)|\leq K(1\wedge
|x|^{2})\quad\Longrightarrow\quad\mathbf{E}(|g(Z_{\Delta}(\alpha)|)\leq
CK\zeta\Delta^{\frac{\beta-\alpha}{\beta}}. \label{eq:Z'sym}%
\end{equation}

\end{lemma}

\begin{proof}
It is similar to the proof of Lemma 2 of \cite{yacjacod04a}. Taking $\eta>0$,
we set $Y_{t}^{\prime\prime}=\sum_{s\leq t}\Delta Y_{s}1_{\{|\Delta
Y_{s}|>\eta\}}$ and $Y^{\prime}=Y-Y^{\prime\prime}$ and if $G\in
\overline{\mathcal{G}}^{\prime}(\zeta,\alpha)$ then $Y$ is symmetrical and
thus we have (47) of the afore--mentioned proof (with $\phi_{\alpha}$
substituted with a constant proportional to $\zeta$), that is
\[
\mathbf{E}(|Y_{\Delta}^{\prime}|^{2})\leq C\zeta\Delta\eta^{2-\alpha}%
\]
for a constant $C$ depending on $\alpha$ only. We also have $Z_{\Delta}%
(\alpha)=\Delta^{-1/\beta}Y_{\Delta}$, hence $|g(Z_{\delta}(\alpha))|\leq
K\Delta^{-2/\beta}|Y_{\Delta}^{\prime}|^{2}$ on the set $\{Y_{\Delta}%
^{\prime\prime}=0\}$, whose probability is smaller than $C\zeta\Delta
/\eta^{\alpha}$. Since $|g|\leq K$, we deduce
\[
\mathbf{E}(|g(Z_{\Delta}(\alpha))|)\leq CK\zeta\left(  \Delta\eta^{-\alpha
}+\Delta^{1-2/\beta}\eta^{2-\alpha}\right)  .
\]
Then take $\eta=\Delta^{1/\beta}$ to obtain the result.
\end{proof}

Next, as soon as the function $k$ satisfies the first half of (\ref{eq:k})
with some $\gamma\geq0$ which has $\gamma<\beta$ whenever $\beta<2$, we set
for $u>0$ and $z\in R$:
\begin{equation}
\Psi_{k}(u,z)=\int h_{\beta}(x)k(ux+z)~dx=\frac{1}{u}\int h_{\beta}\left(
\frac{x}{u}\right)  ~k(x+z)~dx=\frac{1}{u}\int h_{\beta}\left(  \frac{x-z}%
{u}\right)  ~k(x)~dx. \label{eq:P1}%
\end{equation}
(so $\Psi_{k}(u,z)=\Psi_{G,\Delta,\alpha,k}(u,0,z)$, which depends neither on
$G$, nor on $\Delta$, nor on $\alpha$).

\begin{lemma}
\label{lem:2} a) Let $k$ satisfy the first half of (\ref{eq:k}) with some
$\gamma\geq0$ which has $\gamma<\beta$ whenever $\beta<2$. Then $\Psi_{k}$ is
$C^{\infty}$ on $(0,\infty)\times R$. If further $\gamma>0$ and $\nu
\in(0,\infty)$ and $k_{\nu}(x)=k(x)1_{\{|k(x)|\leq\nu\}}$, then for all $K>0$
there exists $M_{K,k}$ such that
\begin{align}
&  |z|\leq K,~~\nu\geq M_{K,k}\qquad\Longrightarrow\qquad\label{eq:FF25}\\
&  \left\vert \frac{\partial^{j+l}}{\partial u^{j}~\partial z^{l}}~\Psi
_{k}(u,z)-\frac{\partial^{j+l}}{\partial u^{j}~\partial z^{l}}~\Psi_{k_{\nu}%
}(u,z)\right\vert \leq\left\{
\begin{array}
[c]{ll}%
C_{j,l,k,K}~u^{\beta-j}~\nu^{1-(l+\beta)/\gamma}\qquad & \mbox{if
}~\beta<2\\
C_{j,l,k,K}~u^{j+l+\gamma-1}~e^{-\nu^{1/\gamma}/u}\qquad & \mbox{if
}~\beta=2.
\end{array}
\right. \nonumber
\end{align}

b) If $k$ is bounded, then for all $\eta\in(0,1)$ we have
\begin{equation}
\eta\leq u\leq1/\eta\quad\Longrightarrow\quad\left\vert \frac{\partial^{j+l}%
}{\partial u^{j}~\partial z^{l}}~\Psi_{k}(u,z)\right\vert \leq C_{l,j,\eta
}~\Vert k\Vert. \label{eq:P3}%
\end{equation}

\end{lemma}

\begin{proof}
(a) If $l\in N$, the $j$th derivative of $u\mapsto(-1)^{l}h_{\beta}%
^{(l)}(x/u)/u^{l+1}$ takes the form $h_{l,j}(x/u)/u^{j+l+1}$ for a function
$h_{l,j}$ satisfying
\begin{equation}
|h_{l,j}(x)|\leq\left\{
\begin{array}
[c]{ll}%
C_{j,l}/(1+|x|^{1+l+\beta})\qquad & \mbox{if }~\beta<2\\
C_{j,l}(1+|x|^{2j+2l})~e^{-x^{2}/2}\qquad & \mbox{if }~\beta=2.
\end{array}
\right.  \label{eq:H}%
\end{equation}
In particular the estimate for $\beta<2$ above also holds for $\beta=2$, and
further $h_{l,j}$ is differentiable and, for all $\beta\in(0,2]$,
\begin{equation}
|h_{l,j}^{\prime}(x)|\leq\frac{C_{j,l}}{1+|x|^{2+l+\beta}}. \label{eq:H'}%
\end{equation}
Therefore we easily deduce from (\ref{eq:P1}) that $\Psi_{k}$ is $C^{\infty}$,
with (by differentiating $l$ times the last term in (\ref{eq:P1}), then $j$
times the analogue of the third term with $h_{\beta}^{(l)}$ instead of
$h_{\beta}$):
\begin{equation}
\frac{\partial^{j+l}}{\partial u^{j}~\partial z^{l}}~\Psi_{k}(u,z)=\frac
{1}{u^{j+l+1}}\int h_{l,j}(x/u)~k(x+z)~dx=\frac{1}{u^{j+l}}\int h_{l,j}%
(x)~k(ux+z)~dx. \label{eq:P2}%
\end{equation}

In particular, for some $\varepsilon_{k}>0$ depending on the function $k$, we
have
\begin{align}
\left\vert \frac{\partial^{j+l}}{\partial u^{j}~\partial z^{l}}~\Psi
_{k}(u,z)-\frac{\partial^{j+l}}{\partial u^{j}~\partial z^{l}}~\Psi_{k_{\nu}%
}(u,z)\right\vert  &  \leq\frac{1}{u^{j+l+1}}\int\big|k(x+z)-k_{\nu
}(x+z)\big|h_{l,j}(x/u)~dx\nonumber\\
&  \leq\frac{C_{k}}{u^{j+l+1}}~\int_{\{1+|x+z|^{\gamma}>\nu\varepsilon_{k}%
\}}(1+|x+z|^{\gamma})~h_{l,j}(x/u)~dx.\nonumber
\end{align}
Then a simple computation, using (\ref{eq:H}), gives us (\ref{eq:FF25}).

(b) When $k$ is bounded, (\ref{eq:H}) and (\ref{eq:P2}) immediately yield
(\ref{eq:P3}).
\end{proof}

Finally we give estimates for the difference $\Psi_{G,\Delta,\alpha,k}$ and
$\Psi_{k}$.

\begin{lemma}
\label{lem:3}If $k$is a bounded function, $\Psi_{G,\Delta,\alpha,k}(u,v,z)$ is
$C^{\infty}$ in $(u,z)$, and for any $\eta\in(0,1)$ we have
\begin{equation}
\eta\leq u\leq1/\eta\quad\Longrightarrow\quad\left\vert \frac{\partial^{j+l}%
}{\partial u^{j}~\partial z^{l}}~\Psi_{G,\Delta,\alpha,k}(u,v,z)\right\vert
\leq C_{l,j,\eta}~\frac{\Vert k\Vert}{1+|z|^{l+\beta}}. \label{eq:P4}%
\end{equation}
Moreover, for all $\eta\in(0,1)$ we have the following, for all $\Delta\leq1$
and $z\in R$ and $u\in\lbrack\eta,1/\eta]$ and $v\in(0,1/\eta]$:

(i) If $G\in G(\phi,\alpha)$ (resp. $G\in\overline{\mathcal{G}}(\zeta,\alpha
)$), then with $\phi_{\alpha}$ given by (\ref{eq:phialp}) (resp. $\phi
_{\alpha}\equiv\zeta$):
\begin{equation}
\left\vert \frac{\partial^{j}}{\partial u^{j}}~\Psi_{G,\Delta,\alpha
,k}(u,v,z)-\frac{\partial^{j}}{\partial u^{j}}~\Psi_{k}(u,0)\right\vert \leq
C_{j,\eta}~\Vert k\Vert~\left(  |z|+\Delta^{\frac{2(\beta-\alpha)}%
{\beta(2+\alpha)}}\phi_{\alpha}(\Delta^{\frac{2+\beta}{\beta(2+\alpha)}%
})\right)  , \label{FF13}%
\end{equation}

(ii) If $G\in\overline{\mathcal{G}}^{\prime}(\zeta,\alpha)$, then
\begin{equation}
\left\vert \frac{\partial^{j}}{\partial u^{j}}~\Psi_{G,\Delta,\alpha
,k}(u,v,z)-\frac{\partial^{j}}{\partial u^{j}}~\Psi_{k}(u,0)\right\vert \leq
C_{j,\eta}~\Vert k\Vert~\left(  |z|+\zeta\Delta^{\frac{\beta-\alpha}{\beta}%
}\right)  , \label{FF13sym}%
\end{equation}

\end{lemma}

\begin{proof}
Observe that $\Psi_{G,\Delta,\alpha,k}(u,v,z)=\int G_{\Delta,\alpha}^{\prime
}(dw)~\Psi_{k}(u,vw+z)$. Then by (\ref{eq:P3}), $\Psi_{G,\Delta,\alpha,k}$ is
$C^{\infty}$ in $(u,z)$, with
\begin{equation}
\frac{\partial^{j+l}}{\partial u^{j}~\partial z^{l}}~\Psi_{G,\Delta,\alpha
,k}(u,v,z)=\int G_{\Delta,\alpha}^{\prime}(dw)~\frac{\partial^{j+l}}{\partial
u^{j}~\partial z^{l}}~\Psi_{k}(u,vw+z), \label{FF11}%
\end{equation}
and for any $\eta\in(0,1)$ we have (\ref{eq:P4}).

Next we prove (i). (\ref{eq:H'}) yields
\begin{equation}
|y|\leq1\quad\Longrightarrow\quad|h_{0,j}(x+y)-h_{0,j}(x)|\leq C_{j,m}%
~\frac{|y|}{1+|x|^{2+\beta}}. \label{eq:zeta}%
\end{equation}
Recalling (\ref{eq:P2}) and (\ref{FF11}), we have
\begin{equation}
\frac{\partial^{j}}{\partial u^{j}}~\Psi_{G,\Delta,\alpha,k}(u,v,z)-\frac
{\partial^{j}}{\partial u^{j}}~\Psi_{k}(u,z)=\int G_{\Delta,\alpha}^{\prime
}(dw)~g(w), \label{eq:dif}%
\end{equation}
where
\begin{align*}
g(w)  &  =\frac{\partial^{j}}{\partial u^{j}}~\Psi_{k}(u,vw+z)-\frac
{\partial^{j}}{\partial u^{j}}~\Psi_{k}(u,z)\\
&  =\frac{1}{u^{j}}\int h_{0,j}(x)\left(  k(ux+vw+z)-k(ux+z)\right)  ~dx\\
&  =\frac{1}{u^{j}}\int\left(  h_{0,j}\left(  x-\frac{vw}{u}\right)
-h_{0,j}(x)\right)  k(ux+z)~dx~,
\end{align*}
for $u,v,z,j$ fixed. Let $\eta\in(0,1)$, and suppose that $\eta\leq
u\leq1/\eta$ and that $v\leq1/\eta$. If $|w|\leq1$ (\ref{eq:zeta}) obviously
yields $|g(w)|\leq C_{j,\eta}~\Vert k\Vert~|w|$, whereas (\ref{eq:P3}) yields
$|g(w)|\leq C_{j,\eta}~\Vert k\Vert$ always: so we have $|g(w)|\leq C_{j,\eta
}~\Vert k\Vert(|w|\bigwedge1)$, and in view of (\ref{eq:dif}) we readily
deduce from (\ref{eq:Z'}) if $G\in G(\phi,\alpha)$ and (\ref{eq:Z'modif}) if
$G\in\overline{\mathcal{G}}(\zeta,\alpha)$ (then $\phi_{\alpha}\equiv\zeta$),
then
\begin{equation}
\left\vert \frac{\partial^{j}}{\partial u^{j}}~\Psi_{G,\Delta,\alpha
,k}(u,v,z)-\frac{\partial^{j}}{\partial u^{j}}~\Psi_{k}(u,z)\right\vert \leq
C_{j,\eta}~\Vert k\Vert~\Delta^{\frac{2(\beta-\alpha)}{\beta(2+\alpha)}}%
\phi_{\alpha}(\Delta^{\frac{2+\beta}{\beta(2+\alpha)}}). \label{eq:P5}%
\end{equation}
Moreover (\ref{eq:P3}) yields $\left\vert \frac{\partial^{j}}{\partial u^{j}%
}~\Psi_{k}(u,z)-\frac{\partial^{j}}{\partial u^{j}}~\Psi_{k}(u,0)\right\vert
\leq C_{j,\eta}~\Vert k\Vert~|z|$, so putting all these together gives
(\ref{FF13}).

Finally we prove (ii). The function $h_{0,j}$ is $C^{\infty}$ and all its
derivatives satisfy the estimates (\ref{eq:H}), and in particular
$H(x)=\sup_{y\in\lbrack x-1/\eta^{2},x+1/\eta^{2}]}|h_{0,j}^{\prime\prime
}(y)|$ is integrable, as well as $h_{0,j}^{\prime}$. Now we have
\begin{equation}
|w|\leq1\quad\Rightarrow\quad\left\vert h_{0,j}\left(  x-\frac{vw}{u}\right)
-h_{0,j}(x)-h_{0,j}^{\prime}(x)\frac{vw}{u}\right\vert \leq C_{j,\eta}%
w^{2}H(x) \label{eq:dif1}%
\end{equation}
as soon as $v<1/\eta$ and $\eta\leq u\leq1/\eta$. Therefore we can write
$g=g_{1}+g_{2}$, where
\[
g_{1}(w)=\frac{vw}{u^{j+1}}~1_{\{|w|\leq1\}}\int h_{0,j}^{\prime
}(x)k(ux+z)dx,
\]%
\[
g_{2}(w)=g(w)1_{\{|w|>1\}}++1_{\{|w|\leq1\}}\int\left(  h_{0,j}\left(
x-\frac{vw}{u}\right)  -h_{0,j}(x)-h_{0,j}^{\prime}(x)\frac{vw}{u}\right)
k(ux+z)dx.
\]
On the one hand, if $G\in\overline{\mathcal{G}}^{\prime}(\zeta,\alpha)$ then
$G_{\Delta,\alpha}^{\prime}$ is symmetrical about $0$, hence $\int
g_{1}(w)G_{\Delta,\alpha}^{\prime}(dw)=0$ because $g_{2}$ is bounded and odd.
On the other hand, (\ref{eq:dif1}) plus the integrability of $H$ and the fact
that $|g(w)|\leq C_{j,\eta}\Vert k\Vert$ yield $|g_{2}(w)|\leq C_{j,\eta}\Vert
k\Vert(w^{2}\bigwedge1)$. Hence, using Lemma \ref{lem:1} we get instead of
(\ref{eq:P5}) that
\begin{equation}
\left\vert \frac{\partial^{j}}{\partial u^{j}}~\Psi_{G,\Delta,\alpha
,k}(u,v,z)-\frac{\partial^{j}}{\partial u^{j}}~\Psi_{k}(u,z)\right\vert \leq
C_{j,\eta}~\Vert k\Vert~\zeta\Delta^{\frac{\beta-\alpha}{\beta}},
\label{eq:P5sym}%
\end{equation}
and we conclude (\ref{FF13sym}) as previously.
\end{proof}

\subsection{Proof of Theorem \ref{theorem:Gknown}.\label{sspr}}

We start by proving (b). With the notation $H=\breve{h}_{\beta}/h_{\beta}$, we
observe that in addition to (\ref{eq:FF}), we have
\[
I(k)=\mathbf{E}(k(W_{1})H(W_{1})),\quad\mathcal{I}(\beta)=\mathbf{E}%
(H(W_{1})^{2}).
\]
An integration by parts yields $E(H(W_{1}))=0$, so $J(k)=E(k^{\prime}%
(W_{1})^{2})$ and $I(k)=E(k^{\prime}(W_{1})H(W_{1}))$ if $k^{\prime
}(x)=k(x)-E(k(W_{1}))$. The desired inequality, which is $I(k)^{2}\leq
J(k)I(\beta)$, follows from the Cauchy--Schwarz inequality. If $k=\overline
{h}_{\beta}$ we also have $k=1+H$, so this inequality is obviously an equality.

For (a), and since $p_{n}\sim n$, we apply Theorem \ref{TTT1}--(b) with
$\chi_{i}^{n}$ given by (\ref{eq:Q1}) and thus $P_{n,\sigma}=P_{\sigma,G}$.
The first step consists in proving (A1) for $S_{n}=S_{n}(G)$. This amounts to
the following lemma, where $\sigma_{n}\rightarrow\sigma>0$ and $P_{n}%
=P_{\sigma_{n},G}$:

\begin{lemma}
\label{LL4} The sequence $S_{n}$ converges to $\sigma$ in probability.
\end{lemma}

\begin{proof}
By (\ref{eq:Z'}) the variables $Z_{\Delta_{n}}^{n}(\beta)$ associated with the
law $G^{n}$ converge in law to $0$ (because $\phi_{\beta}(x)\rightarrow$ as
$x\rightarrow0$). The variables $\chi_{i}^{\prime n}$, which equal $\sigma
_{n}W_{1}+Z_{\Delta_{n}}^{n}(\beta)$ in law, converge in law to $\sigma W_{1}%
$. Hence $\gamma_{n}:=P_{n}(|\chi_{i}^{\prime n}|>1)\rightarrow\psi(\sigma)$.
If $\zeta_{i}^{n}=1_{\{|\chi_{i}^{\prime n}|>1\}}$, (\ref{LLN}) applied with
$q_{n}=m_{n}$ yields $V_{n}~\overset{\mathbf{P}_{\sigma,G}}{\longrightarrow
}~\psi(\sigma)$. Since $\psi^{-1}$ is $C^{\infty}$ and strictly monotone, the
result readily follows.
\end{proof}

Next we set $Q_{n}=0$, so (A2) is satisfied, and
\[
f_{n,s,q}(x)=k_{n}\left(  \frac{\Delta_{n}^{-1/\beta}(x-b^{\prime}%
(G,\beta)\Delta_{n}))}{s}\right)  ,\qquad H_{n,s}(u)=\Psi_{G,\Delta_{n}%
,\beta,k_{n}}\left(  \frac{u}{s},\frac{1}{s},0\right)  .
\]
Upon comparing (\ref{FE6}) and (\ref{FE7}) with (\ref{Q12}) and (\ref{Q13}),
we see that $\widehat{\sigma}_{n}(G,\phi,k)=\widehat{\sigma}_{n}(S_{n},Q_{n}%
)$. Therefore it remains to prove (B1)--(B6) with a sequence $w_{n}$
satisfying $w_{n}/\sqrt{p_{n}}\rightarrow\infty$, and that
\begin{equation}
\Xi^{2}\sigma)=\sigma^{2}J(k)/I(k)^{2}. \label{eq:Q25}%
\end{equation}

Observe that under $P_{\sigma,G}$ the variables $\chi_{i}^{n}$ have the same
law as $\sigma W_{1}+Z_{\Delta_{n}}(\beta)$. Then (\ref{Q11}) gives
$F_{n,s,q}(\sigma)=H_{n,s}(\sigma)$. It follows that (B6) holds with $w_{n}$
arbitrarily large, while (B2) follows from (\ref{FF11}).

If $k$ is bounded, hence $k_{n}=k$, we have $\Vert f_{n,s}\Vert\leq\Vert
k\Vert$ and (B1) is obvious; further, (\ref{FF13}) with $\alpha=\beta$ and
$k^{r}$ yields
\[%
\begin{array}
[c]{l}%
j=0,1,~~r=1,2,~~\eta\leq u\leq\frac{1}{\eta},~~v\leq\frac{1}{\eta}%
\qquad\Longrightarrow\\
\left\vert \frac{\partial^{j}}{\partial u^{j}}~\Psi_{G,\Delta_{n},\beta,k^{r}%
}(u,v,0)-\frac{\partial^{j}}{\partial u^{j}}~\Psi_{k^{r}}(u,0)\right\vert \leq
C_{\eta,k}\phi_{\beta}(\Delta_{n}^{1/\beta}),
\end{array}
\]
which gives (B3) with $\overline{F}_{s}(u)=\Psi_{k}(u/s,0)$ and (B5) with
$F^{(2)}(u)=\Psi_{k^{2}}(1,0)$. On the other hand when $k$ is unbounded we
have $\Vert f_{n,s}\Vert\leq\nu_{n}$ and thus (B1) follows from (\ref{eq:kn});
further, $\nu_{n}\rightarrow\infty$ and we can combine (\ref{FF13}) with
(\ref{eq:FF25}) to get for all $n$ large enough:
\[%
\begin{array}
[c]{l}%
j=0,1,~~r=1,2,~~\eta\leq u\leq\frac{1}{\eta},~~v\leq\frac{1}{\eta}%
\qquad\Longrightarrow\\
\left\vert \frac{\partial^{j}}{\partial u^{j}}~\Psi_{G,\Delta_{n},\beta
,k_{n}^{r}}(u,v,0)-\frac{\partial^{j}}{\partial u^{j}}~\Psi_{k^{r}%
}(u,0)\right\vert \leq\left\{
\begin{array}
[c]{ll}%
C_{\eta,k}\left(  \nu_{n}^{r}\phi_{\beta}(\Delta_{n}^{1/\beta})+\frac{1}%
{\nu_{n}^{\beta/r\gamma-1}}\right)  \qquad & \mbox{if }~\beta<2\\
C_{\eta,k}\left(  \nu_{n}^{r}\phi_{2}(\Delta_{n}^{1/2})+e^{-\eta\nu
_{n}^{1/r\gamma}}\right)  \qquad & \mbox{if }~\beta=2.
\end{array}
\right.
\end{array}
\]
Then, in view of (\ref{eq:kn}) and $2\gamma<\beta$ when $\beta<2$, we again
deduce (B3) with $F_{s}(u)=\Psi_{k}(u/s,0)$ and (B5) with $F^{(2)}%
(u)=\Psi_{k^{2}}(1,0)$.

Since $h_{0,1}=-\breve{h}_{\beta}$, we deduce that $\overline{F}_{\sigma
}^{\prime}(\sigma)=\Psi_{k}^{\prime}(1,0)/\sigma=-I(k)/\sigma$ (recall
(\ref{eq:P2}) and the second part of (\ref{eq:k})), hence (B4) holds. We also
have $\overline{F}_{\sigma}(\sigma)=\Psi_{k}(1,0)=E(k(W_{1}))$ and
$F^{(2)}(\sigma)=E(k(W_{1})^{2})$, hence $J(k)=F^{(2)}(\sigma)-F_{\sigma
}(\sigma)^{2}$ and (\ref{eq:Q25}) follows.

\section{Estimation of $\sigma$ in the semiparametric
case.\label{sec:SemiParEst}}

Perhaps more realistic than the situation of Theorem \ref{theorem:Gknown} is
the case where we want to estimate $\sigma,$ but the measure $G$ is unknown,
although we know that it belongs to the class $G_{\beta}$. This is a
semiparametric situation: parametric as far as $\sigma W_{t}$ is concerned,
but nonparametric as far as $Y_{t}$ is concerned. Because $G$ is unknown, the
estimating equations in this case must be based on the law of $W$ alone. The
challenge is then to achieve rate efficiency despite the lack of information
about $G.$

\subsection{Construction of the estimators.}

As said before, we cannot hope for estimators $\widehat{\sigma}_{n}$ that
behave nicely for all $G\in G_{\beta}$ at once. Therefore we suppose that $G$
is unknown, but is known to belong to $\overline{\mathcal{G}}(\zeta,\alpha)$
for some $\alpha<\beta$ and some $\zeta>0$: we refer to this as Case 1. We
also consider a more restrictive situation, called Case 2, for which $G$ is
known to belong to the set $\overline{\mathcal{G}}^{\prime}(\zeta,\alpha)$.

The construction looks pretty much like the previous one, except that besides
our preliminary estimator for $\sigma$ we need to produce an estimator $B_{n}$
for the drift $b^{\prime}(G,\alpha)$ in order to remove it. In Case 2, since
we know that $b^{\prime}(G,\alpha)=0$ we just set
\begin{equation}
B_{n}=0. \label{FE2'}%
\end{equation}
In Case 1 we set $m_{n}=[\delta n]$ for some arbitrary $\delta\in(0,1/2)$
($[x]$ denotes the integer part of $x$), so that $m_{n}\sim\delta n$. Then we
pick a $C^{\infty}$ and strictly increasing and odd function $\theta$, with
bounded derivative and $\theta(0)=0$ and $\theta(\pm\infty)=\pm1$ (for example
$\theta(x)=\frac{2}{\pi}\arctan(x)$~), and set for $u\in R$
\begin{equation}
R_{n}(u)=\frac{1}{m_{n}}\sum_{i=1}^{m_{n}}\theta(\Delta_{n}^{-1/\beta}%
(\chi_{i}^{n}-u)). \label{NP21}%
\end{equation}
Since $u\mapsto R_{n}(u)$ is continuous and decreases strictly from $+1$ to
$-1$ as $u$ goes from $-\infty$ to $+\infty$, we can set
\begin{equation}
B_{n}=\inf(u:R_{n}(u)=0)\quad(=~\mbox{the only root
of $R_n(.)=0$}~). \label{NP22}%
\end{equation}

Next we construct our preliminary estimator for $\sigma$. In Case 1, and with
$m_{n}$ as above, we set $q_{n}=m_{n}$ and $p_{n}=n-2m_{n}$. In Case 2, we
choose a sequence $m_{n}$ satisfying (\ref{FE1}) and then we set $q_{n}=0$ and
$p_{n}=n-m_{n}$. Then in both cases we set
\begin{equation}
V_{n}=\frac{1}{m_{n}}~\sum_{i=q_{n}+1}^{q_{n}+m_{n}}1_{\{|\Delta_{n}%
^{-1/\beta}(\chi_{i}^{n}-B_{n})|>1\}}\label{NP23}%
\end{equation}
and%
\begin{equation}
S_{n}=\left\{
\begin{array}
[c]{ll}%
\psi^{-1}(V_{n})\qquad & \mbox{if }~0<V_{n}<1\\[2mm]%
1 & \mbox{otherwise}.
\end{array}
\right.  \label{NP24}%
\end{equation}

To form estimating equations for $\sigma$, we choose a function $k$ satisfying
(\ref{eq:k}) with $\gamma=0$ (that is, $k$ is bounded and $I(k)\neq0$). With
$\Psi_{k}$ given by (\ref{eq:P1}) we define the estimating functions (for
$u>0$)
\begin{equation}
U_{n}(u)=\frac{1}{p_{n}}\sum_{i=q_{n}+m_{n}+1}^{n}k\left(  \frac{\Delta
_{n}^{-1/\beta}(\chi_{i}^{n}-B_{n})}{S_{n}}\right)  -\Psi_{k}\left(  \frac
{u}{S_{n}},0\right)  ,\label{FE6'}%
\end{equation}
and the final estimators
\begin{equation}
\widehat{\sigma}_{n}(k)=\left\{
\begin{array}
[c]{ll}%
\mbox{the $u$ with $U_n(u)=0$ which is closest to ~$S_n$}\quad &
\mbox{if it exists}\\[2.5mm]%
1 & \mbox{otherwise}.
\end{array}
\right.  \label{FE7'}%
\end{equation}

Note that, unlike the centering $\Psi_{G,\Delta_{n},\beta,k_{n}}\left(
\frac{u}{S_{n}(G)},\frac{1}{S_{n}(G)},0\right)  $ utilized in the parametric
case (recall (\ref{FE6})), the centering we now use, based on $\Psi_{k}\left(
\frac{u}{S_{n}},0\right)  $ in (\ref{FE6'}) does not involve the measure $G.$
Indeed, these estimators depend explicitly on $\beta$ and $k$, but on nothing
else, and in particular not on $G$. Observe that they are much easier to
compute than the estimator of the parametric case. This is particularly true
when $k(x)=\cos(wx)$ for some $w>0$, since then $\Psi_{k}(u,0)=e^{-w^{\beta
}u^{\beta}/2}$ is invertible in $u$, and we will detail this example in the
next section, but it is also true in general: first because they depend only
on the function $\Psi_{k}(u,.)$ which is much simpler than the function
$\Psi_{G,\Delta,\beta,k}$ accruing in the estimation in the parametric case,
second because as a rule $u\mapsto\Psi_{k}(u,0)$ is at least \textquotedblleft
locally invertible\textquotedblright\ around $u=1$.

The estimators (\ref{FE6'}) have formally the same expression in both Case 1
and Case 2, but the preliminary estimators $B_{n}$ and $S_{n}$ disagree for
the two cases and also $p_{n}\sim(1-2\delta)n$ in Case 1 and $p_{n}\sim n$ in
Case 2, a difference which is important for the asymptotic variance of the
estimators. So we will write \textquotedblleft the Case 1
version\textquotedblright\ or \textquotedblleft the Case 2
version\textquotedblright\ of the estimator.

\subsection{Asymptotic distribution in the semiparametric case.}

Recall the notation $I(k)$ and $J(k)$ and $\Sigma^{2}(k)$ of (\ref{eq:k}) and
(\ref{eq:FF}), and let us add some other:
\begin{equation}
\rho(\alpha,\beta)=\frac{2(\beta-\alpha)}{\beta(2+\alpha)},\qquad\rho^{\prime
}(\alpha,\beta)=\frac{\beta-\alpha}{\beta}.\label{eq:rho}%
\end{equation}
Observe that $\rho(\alpha,\beta)<\rho^{\prime}(\alpha,\beta)$ always.

\begin{theorem}
\label{theorem:Gunknown1} Let $\alpha\in(0,\beta)$ and $\zeta>0$, and $k$ be a
bounded function with $I(k)\neq0$, and $\varepsilon\in(0,1)$. Take the Case 1
version of the estimators.

a) If
\begin{equation}
\sup_{n}~n\Delta_{n}^{2\rho(\alpha,\beta)}\rightarrow0, \label{eq:aseff}%
\end{equation}
the sequence $\sqrt{n}~(\widehat{\sigma}_{n}(k)-\sigma)$ converges in law to
$N(0,\sigma^{2}\Sigma^{2}(k)/(1-2\delta))$ under $P_{\sigma,G}$, uniformly in
$n\geq1$ and in $\sigma\in\lbrack\varepsilon,1/\varepsilon]$ and in
$G\in\overline{\mathcal{G}}(\zeta,\alpha)$.

b) In general, the variables $(\sqrt{n}\bigwedge\Delta_{n}^{-\rho(\alpha
,\beta)})(\widehat{\sigma}_{n}(k)-\sigma)$ are tight under $P_{\sigma,G}$,
uniformly in $\sigma\in\lbrack\varepsilon,1/\varepsilon]$ and in
$G\in\overline{\mathcal{G}}(\zeta,\alpha)$ and $n$.
\end{theorem}

\begin{theorem}
\label{theorem:Gunknown} Let $\alpha\in(0,\beta)$ and $\zeta>0$, and $k$ be a
bounded function with $I(k)\neq0$, and $\varepsilon\in(0,1)$. Take the Case 2
version of the estimators.

a) If
\begin{equation}
\sup_{n}~n\Delta_{n}^{2\rho^{\prime}(\alpha,\beta)}\rightarrow0,
\label{eq:aseffsym}%
\end{equation}
the sequence $\sqrt{n}~(\widehat{\sigma}_{n}(k)-\sigma)$ converges in law to
$N(0,\sigma^{2}\Sigma^{2}(k))$ under $P_{\sigma,G}$, uniformly in $n\geq1$ and
in $\sigma\in\lbrack\varepsilon,1/\varepsilon]$ and in $G\in\overline
{\mathcal{G}}^{\prime}(\zeta,\alpha)$.

b) In general, the variables $(\sqrt{n}\bigwedge\Delta_{n}^{-\rho^{\prime
}(\alpha,\beta)})(\widehat{\sigma}_{n}(k)-\sigma)$ are tight under
$P_{\sigma,G}$, uniformly in $\sigma\in\lbrack\varepsilon,1/\varepsilon]$ and
in $G\in\overline{\mathcal{G}}^{\prime}(\zeta,\alpha)$ and $n$.
\end{theorem}

The optimal\ choice of the function $k$ has been discussed after Theorem
\ref{theorem:Gknown}: when $\beta<2$, we have asymptotic efficiency in the
situation of the second theorem above, provided we take $k=\overline{h}%
_{\beta}$, and despite the fact that we are in a semiparametric setting. When
$\beta=2$ the choice $k=\overline{h}_{\beta}$, that is $k(x)=-x^{2}$, is not
permitted in the above theorem, but with $k(x)=-x^{2}1_{\{|x|\leq A\}}$ one
achieves an asymptotic variance which approaches the optimal variance when $A$
goes to infinity: see Section \ref{sec:Example2}.

Also, some other comments are in order here:

\begin{remark}
When $\alpha$ increases, then $\rho(\alpha,\beta)$ and $\rho^{\prime}%
(\alpha,\beta)$ decrease, so (\ref{eq:aseff}) and (\ref{eq:aseffsym}) are more
difficult to obtain and the \textquotedblleft rate\textquotedblright\ in (b)
of the two theorems above gets worse, as it should be.
\end{remark}

\begin{remark}
In connection with what precedes, one should mention that when (\ref{eq:aseff}%
) fails the actual rate of convergence (that is, a sequence $\delta_{n}$ such
that the law of $\delta_{n}((\widehat{\sigma}_{n}(k)-\sigma)$ converges to a
non--degenerate limit, or at least admits among its weak limiting measures a
non--degenerate one) is not only unknown, but actually depends on the true
underlying (unknown) measure $G$ and in particular on the minimal index
$\alpha^{\prime}$ such that $G\in\overline{\mathcal{G}}_{\alpha^{\prime}}$ (we
know that $\alpha^{\prime}\leq\alpha$, but the inequality could be strict). In
other words, the rate could be for example $\sqrt{n}$ for a particular $G$,
even without (\ref{eq:aseff}).
\end{remark}

\begin{remark}
However we will see in the examples below (see Section \ref{sec:Example4} in
particular) that (\ref{eq:aseffsym}) is necessary for having convergence to a
centered distribution with rate $\sqrt{n}$ and also that the rate in (b) of
Theorem \ref{theorem:Gunknown} is sharp, if we want to have a result which
holds uniformly in $G\in\overline{\mathcal{G}}^{\prime}(\zeta,\alpha)$. We do
not know whether (\ref{eq:aseff}) or the rate in (b) are optimal for Theorem
\ref{theorem:Gunknown1}.
\end{remark}

\begin{remark}
Of course it might exist other -- thoroughly different -- estimators behaving
better than the $\widehat{\sigma}_{n}(k)$'s, and perhaps having a better rate
than in (b) of these theorems (the rate cannot be improved in (a), of course).
We think this doubtful, however.
\end{remark}

\begin{remark}
The most interesting situation is when we have asymptotic efficiency (this
happens when $G$ is symmetrical), or at least \textquotedblleft
rate--efficiency\textquotedblright\ (that is of order $\sqrt{n}$). We have
this under (\ref{eq:aseff}) or (\ref{eq:aseffsym}), which mean that
$\Delta_{n}$ goes to $0$ fast enough. Of course having $\Delta_{n}=o(1/n)$ is
of no practical use. When $\Delta_{n}=1/n$, then rate--efficiency is satisfied
as soon as $\alpha\leq2\beta/(4+\beta)$ for the first theorem and $\alpha
\leq\beta/2$ for the second one. If $Y$ is a compound Poisson process with
drift, rate efficiency holds as soon as $n\Delta_{n}^{2}$ is bounded, whatever
$\beta\in(0,2]$ is (take $\alpha=0$).
\end{remark}

\begin{remark}
When we do not know that $G$ is symmetrical we cannot achieve asymptotic
efficiency even under (\ref{eq:aseff}). However the asymptotic variances in
the two theorems above are the same, up to the factor $1-2\delta$: hence by
choosing $\delta$ small one can approach asymptotic efficiency as much as one
wants to.
\end{remark}

\subsection{Proof of Theorems \ref{theorem:Gunknown1} and
\ref{theorem:Gunknown}.\label{sspr'}}

As above, we refer to Theorem \ref{theorem:Gunknown1} as to Case 1, and to
Theorem \ref{theorem:Gunknown} as to Case 2. The proof goes through several steps.

1) We fix $\alpha\in(0,\beta)$ and $\zeta>0$. The sequence $\Delta_{n}$ is
fixed, and we set
\begin{equation}
\rho=\left\{
\begin{array}
[c]{ll}%
\rho(\alpha,\beta)\quad & \mbox{in Case
1}\\[2.5mm]%
\rho^{\prime}(\alpha,\beta)\quad & \mbox{in Case 2},
\end{array}
\right.  \qquad\lambda_{n}=\sqrt{n}\bigwedge\frac{1}{\Delta_{n}^{\rho}}.
\label{eq:rho'}%
\end{equation}
In order to get tightness or convergence, \textquotedblleft
uniform\textquotedblright\ in $\sigma$ and in $G$ is the relevant class, it is
of course enough to take a sequence $\sigma_{n}\rightarrow\sigma>0$ and a
sequence $G^{n}$ in $\overline{\mathcal{G}}(\zeta,\alpha)$ (resp.
$\overline{\mathcal{G}}^{\prime}(\zeta,\alpha)$), and to prove the tightness
or convergence in law of the normalized estimation errors $\widehat{\sigma
}_{n}-\sigma_{n}$, under the measures $P_{n}=P_{\sigma_{n},G^{n}}$. Below we
fix the sequences $\sigma_{n}$ and $G^{n}$.

Finally, we denote by $Z_{n}:=Z_{\Delta_{n}}^{n}(\alpha)$ the variable
associated with the measure $G^{n}$ by (\ref{eq:b'}), and we set
$b_{n}^{\prime}=\Delta_{n}^{1-1/\beta}b^{\prime}(G^{n},\alpha)$, which
vanishes in Case 2. \medskip

2) Let $Q_{n}=\lambda_{n}B_{n}^{\prime}$, where $B_{n}^{\prime}=(\Delta
_{n}^{-1/\beta}B_{n}-b_{n}^{\prime})$. We want to prove that the sequence
$Q_{n}$ satisfies (A2). This is obvious in Case 2 because $Q_{n}=0$. So we
suppose that we are in Case 1. Let us introduce some notation: with $j=1,2$
and $\theta^{\prime}$ being the derivative of $\theta$, we put
\begin{equation}
\Gamma_{j}(\sigma)=\mathbf{E}(\theta(\sigma W_{1})^{j}),\qquad\Gamma
_{1}^{\prime}(\sigma)=\mathbf{E}(\theta^{\prime}(\sigma W_{1})) \label{eq:Gam}%
\end{equation}
($\Gamma_{1}^{\prime}$ is of course the derivative of $\Gamma_{1}$).

Observe that $B_{n}^{\prime}$ is the only root of $R_{n}(.)=0$, where
\[
\mathcal{R}_{n}(u)=R_{n}(\Delta_{n}^{1/\beta}(u+b_{n}^{\prime}))=\frac
{1}{m_{n}}\sum_{i=1}^{m_{n}}\zeta_{i}^{n}(u),\quad\mbox{with }~\zeta_{i}%
^{n}(u)=\theta(\Delta_{n}^{-1/\beta}\chi_{i}^{n}-u-b_{n}^{\prime}).
\]
The $\zeta_{i}^{n}(u)$'s for $i\geq1$ are i.i.d. with the same law (under
$P_{n}$) than the variable $\theta(\sigma_{n}W_{1}+Z_{n}-u)$ (we have used
here the scaling property of $W$).

The functions $\gamma_{n,j}(u)=E_{n}((\zeta_{i}^{n}(u)^{j})$, for $j\in N$,
are $C^{\infty}$ and bounded as well as their derivatives, uniformly in $u$
and $n$, and we can interchange derivation and expectation. So we can apply
(\ref{eq:Z'modif}) to the functions $g_{n,j,p}(w)=\int h_{\beta}%
(x)(\partial^{p}\theta^{j}/\partial u^{p})(\sigma_{n}x+w-u)-(\partial
^{p}\theta^{j}/\partial u^{p})(\sigma_{n}x-u))~dx$, to get for $p,j\in N$:
\begin{equation}
\left\vert \frac{\partial^{p}}{\partial u^{p}}~\gamma_{n,j}(u)-\Gamma
_{j,p}(\sigma_{n},u)\right\vert \leq C_{p,j}\zeta\Delta_{n}^{\rho}%
,\quad\mbox{where }~\Gamma_{j_{p}}(v,u)=(-1)^{p}\int\frac{\partial^{p}%
\theta^{j}}{\partial u^{p}}(vx-u)h_{\beta}(x)dx. \label{NP24'}%
\end{equation}
In particular $\Gamma_{j,0}(\sigma,0)=\Gamma_{j}(\sigma)$ for $j=1,2$ and
$\Gamma_{1,1}(\sigma_{n},0)=\Gamma_{1}^{\prime}(\sigma)$ with the notation
(\ref{eq:Gam}).

Now, $R_{n}$ also is $C^{\infty}$, bounded as well as all its derivatives,
uniformly in $n$, $u$ and $\omega$. So an application of Lemma \ref{LE3} and
the continuity of the functions $\Gamma_{j,p}$ readily yield
\begin{equation}
\frac{\partial^{p}}{\partial u^{p}}~\mathcal{R}_{n}(u)~\rightarrow
~\Gamma_{1,p}(\sigma,u)~~\mbox{locally
uniformly in $u$, in $\mathbf{P}_n$--probability,} \label{NP25}%
\end{equation}%
\begin{equation}
\eta_{n}:=\sqrt{m_{n}}~(\mathcal{R}_{n}(0)-\gamma_{n,1}(0))~\overset
{{\mathcal{L}}(\mathbf{P}_{n})}{\longrightarrow}~\mathcal{N}(0,\Gamma
_{2}(\sigma)-\Gamma_{1}(\sigma)^{2})). \label{NP26}%
\end{equation}

The properties of $\theta$ imply that $u\mapsto\Gamma_{1,0}(\sigma,.)$
decreases strictly and vanishes at $0$; since by construction $R_{n}%
(B_{n}^{\prime})=0$, we deduce from (\ref{NP25}) for $p=0$ that $B_{n}%
^{\prime}~\overset{\mathbf{P}_{n}}{\longrightarrow}~0$. Another application of
(\ref{NP25}) yields that $R_{n}^{\prime}(B_{n}^{\prime\prime})~\overset
{\mathbf{P}_{n}}{\longrightarrow}~\Gamma_{1}^{\prime}(\sigma)$ for any
sequence $B_{n}^{\prime\prime}$ of random variable going to $0$ in $P_{n}%
$--probability. Since $R_{n}(B_{n}^{\prime})=0$ we have
\begin{equation}
\mathcal{R}_{n}^{\prime}(B_{n}^{\prime\prime})~B_{n}^{\prime}=-\mathcal{R}%
_{n}(0)=-\frac{\eta_{n}}{\sqrt{m_{n}}}-\gamma_{n,1}(0) \label{NP29}%
\end{equation}
for some random variable $B_{n}^{\prime\prime}$ satisfying $|B_{n}%
^{\prime\prime}|\leq|B_{n}^{\prime}|$. Moreover $\Gamma_{1,0}(0)=0$, due to
the fact that $\theta$ is odd, hence $|\gamma_{n,1}(0)|\leq C\zeta\Delta
_{n}^{\rho}$ by (\ref{NP24'}). Since $R_{n}^{\prime}(B_{n}^{\prime\prime
})~\overset{\mathbf{P}_{n}}{\longrightarrow}~\Gamma_{1}^{\prime}(\sigma)\neq
0$, we deduce that $Q_{n}=\lambda_{n}B_{n}^{\prime}$ satisfies (A2) from
(\ref{NP26}) (recall $m_{n}\sim\delta n$ here and (\ref{eq:rho'})). \smallskip

3) Now we proceed to proving the consistency of the preliminary estimators
$S_{n}$. In Case 2 the variables $V_{n}$ and $S_{n}$ are the variables
$V_{n}(G^{n})$ and $S_{n}(G^{n})$ of (\ref{FE2}) and (\ref{FE4}) (they do not
depend on $G^{n}$ in fact), so the result follows from Lemma \ref{LL4}. In
Case 1, set
\[
V_{n}(v)=\frac{1}{m_{n}}~\sum_{i=q_{n}+1}^{q_{n}+m_{n}}1_{\{|\Delta
_{n}^{-1/\beta}(\chi_{i}^{n}-v)|>1\}},\qquad\delta_{n}(v)=\mathbf{P}%
_{n}(|\Delta_{n}^{-1/\beta}(\chi_{i}^{n}-v)|>1).
\]
Then (\ref{LLN}) yields
\begin{equation}
V_{n}(v_{n})-\delta_{n}(v_{n})~\overset{\mathbf{P}_{n}}{\longrightarrow}~0.
\label{NP28}%
\end{equation}
However, $\Delta_{n}^{-1/\beta}(\xi_{i}^{n}-v_{n})$ has the same distribution
as $\sigma_{n}W_{1}+Z_{n}+b_{n}^{\prime}-\Delta_{n}^{-1/\beta}v_{n}$, which by
(\ref{eq:Z'}) converges in law to $\sigma W_{1}$ as soon as $b_{n}^{\prime
}-\Delta_{n}^{-1/\beta}v_{n}\rightarrow0$. Since $B_{n}$ and $(V_{n}(v):v\in
R)$ are independent and $B_{n}^{\prime}=\Delta_{n}^{-1/\beta}B_{n}%
-b_{n}^{\prime}~\overset{\mathbf{P}_{n}}{\longrightarrow}~0$ because
$Q_{n}=\lambda_{n}B_{n}^{\prime}$ satisfies (A2) and $\lambda_{n}%
\rightarrow\infty$, we deduce from (\ref{NP28}) that $V_{n}=V_{n}%
(B_{n})~\overset{\mathbf{P}_{n}}{\longrightarrow}~\psi(\sigma)$. Then the
consistency is proved like in the end of Lemma \ref{LL4}. \smallskip

4) At this stage we will apply Theorem \ref{TTT1}, with the variables
$(S_{n},Q_{n})$ as above and the i.i.d. variables $(\chi_{q_{n}+m_{n}+i}%
^{n}:1\leq i\leq p_{n})$. Observe that with the notation (\ref{Q13}) and
(\ref{FE7'}), we have $\widehat{\sigma}_{n}^{\prime}(k)=\widehat{\sigma}%
_{n}(S_{n},Q_{n})$. We have shown (A1) and (A2) in the two previous steps.
Set
\[
f_{n,s,q}(x)=k\left(  \frac{\Delta_{n}^{-1/\beta}x-b_{n}^{\prime}%
-q/\lambda_{n}}{s}\right)  ,\qquad H_{n,s}(u)=\Psi_{k}\left(  \frac{u}%
{s},0\right)  .
\]
Then (\ref{Q11}) gives for $r=1,2$:
\[
F_{n,s,q}(u)=\Psi_{G^{n},\Delta_{n},\alpha,k}\left(  \frac{u}{s},\frac{1}%
{s},-\frac{q}{s\lambda_{n}}\right)  ,\qquad F_{n,s,q}^{(2)}(u)=\Psi
_{G^{n},\Delta_{n},\alpha,k^{2}}\left(  \frac{u}{s},\frac{1}{s},-\frac
{q}{s\lambda_{n}}\right)  .
\]

Let us check (B1)--(B6). Since $k$ is bounded, (B1) is obvious, whereas (B2)
follows from Lemma \ref{lem:2}. Next, if we set $\overline{F}_{s}(u)=\Psi
_{k}(u/s,0)$ and $F^{(2)}(u)=\Psi_{k^{2}}(1,0)$, Lemma \ref{lem:3} yields for
$j=0,1$ and $\eta\in(0,1)$ and $s,u\in\lbrack\eta,1/\eta]$ and $|q|\leq1/\eta
$:
\[
\left\vert \frac{\partial^{j}}{\partial u^{j}}~H_{n,s}(u)-\frac{\partial^{j}%
}{\partial u^{j}}~\overline{F}_{s}(u)\right\vert \leq C_{k,\eta}\zeta
\Delta_{n}^{\rho},
\]%
\[
\left\vert F_{n,s,q}^{(2)}(u)-F_{s}^{(2)}(u)\right\vert \leq C_{k,\eta
}~\left(  \zeta\Delta_{n}^{\rho}+\frac{1}{\lambda_{n}}\right)  ,
\]%
\[
\left\vert F_{n,s,q}(u)-H_{n,s}(u)\right\vert \leq C_{k,\eta}~\left(
\zeta\Delta_{n}^{\rho}+\frac{1}{\lambda_{n}}\right)  .
\]
These give (B3) and (B5), and also (B6) with $w_{n}=\lambda_{n}$. Finally (B4)
holds because $\overline{F}_{s}^{\prime}(s)=\psi_{k}^{\prime}(1,0)/s=-I(k)/s$,
and (\ref{eq:Q25}) holds here as well as in the previous section.

We can thus apply Theorem \ref{TTT1}: the sequence $\lambda_{n}(\widehat
{\sigma}_{n}-\sigma_{n})$ is tight under $P_{n}$ in all cases, and this gives
the two claims (b). Under (\ref{eq:aseff}) or (\ref{eq:aseffsym}) we have
$\lambda_{n}/\sqrt{n}\rightarrow\infty$, hence $\lambda_{n}\sqrt{p_{n}%
}\rightarrow\infty$ as well, so $\sqrt{p_{n}}~(\widehat{\sigma}_{n}-\sigma
_{n})$ converges in law under $P_{n}$ to a centered Gaussian variable with
variance
\[
\Xi^{2}(\sigma)=\frac{F^{(2)}(\sigma)-\overline{F}_{\sigma}(\sigma)^{2}%
}{\overline{F}_{\sigma}(\sigma)^{2}},
\]
which in view of $\overline{F}_{\sigma}(\sigma)^{2}=J(k)/\sigma^{2}$ equals
$\sigma^{2}\Sigma^{2}(k)$: since $p_{n}\sim(1-2\delta)n$ in Case 1 and
$p_{n}\sim n$ in Case 2, we obtain the two claims (a).

\section{Example: The empirical characteristic function.\label{sec:Example1}}

We now turn to specific estimators. To each specification of an admissible
function $k$ (in the sense of satisfying the assumptions of the above
results), corresponds an estimator for $\sigma.$ For instance, one way of
estimating a parameter for i.i.d. variables $X_{j}$ is to use the empirical
characteristic function, that is $\sum_{j\in J}\exp(iwX_{j})$ for some given
$w$ (or several $w$'s at once) and where $J$ is the index set. If the $X_{j}%
$'s are symmetrical, one should in fact look at the real part only, that is
$\sum_{j\in J}\cos(wX_{j})$. Other estimators based on the empirical
characteristic function in related contexts are given by e.g., \cite{press72},
\cite{fenech76}, \cite{feuervergermcdunnough81a}, Chapter 4 in
\cite{zolotarev} and \cite{singleton01}.

In the parametric situation, at stage $n$ the variable $X_{j}$ is
$\chi^{\prime}{}_{j}^{n}(G)$ and $J=\{m_{n}+1,\ldots,n\}$. Those variables are
\textquotedblleft almost\textquotedblright\ symmetrical (the leading term $W$
coming in them is symmetrical). So we consider for any given $w>0$ the
variable
\begin{equation}
V_{n}(w)=\frac{1}{p_{n}}\sum_{i=m_{n}+1}^{n}\cos\left(  \frac{w\chi^{\prime}%
{}_{i}^{n}(G)}{S_{n}(G)}\right)  , \label{eq:charf}%
\end{equation}
where $S_{n}(G)$ is the preliminary estimator. In other words, if we take
$k(x)=\cos(wx)$ (a bounded function, so $k_{n}=k$ in (\ref{eq:knn})), the
estimating function of (\ref{FE6}) is
\begin{equation}
U_{n,G,\beta,k}(u)=V_{n}(w)-\Psi_{G,\Delta_{n},\beta,k}\left(  \frac{u}%
{S_{n}(G)},\frac{1}{S_{n}(G)},0\right)  . \label{eq:char1}%
\end{equation}

Furthermore, this class of functions $k$ is one for which the function
$\Psi_{G,\Delta,\beta,k}$ is explicit, at least when the exponent in the
L\'{e}vy--Khintchine formula for $Y$ is explicitly known. More precisely, let
us write $\rho(u)$ for the exponent in (\ref{eq:Ychar}), and recall that
$E(\exp iuY_{t})=\exp t\rho(u)$. Then obviously when $g(x)=e^{iwx}$ we have
\[
\Psi_{G,\Delta,\beta,g}(u,v,0)=\exp\left(  -\frac{w^{\beta}u^{\beta}}%
{2}+\Delta\rho(wv\Delta^{-1/\beta})-iwvb^{\prime}(G,\alpha)\Delta^{1-1/\alpha
}\right)  .
\]
Taking the real part, and using (\ref{eq:b'}) and the fact that $G\in
G_{\beta}$, we see that for $k(x)=\cos(wx)$ we have
\begin{equation}
\Psi_{G,\Delta,\beta,k}(u,v,0)=e^{A_{\Delta}(u,v)}~\cos(B_{\Delta
}(u,v)),\label{eq:Ex3}%
\end{equation}
where
\begin{equation}
A_{\Delta}(u,v)=-\frac{w^{\beta}u^{\beta}}{2}+\int F(dx)\left(  \cos
(wv\Delta^{1-1/\beta}x)-1\right)  ,\label{eq:Ex4}%
\end{equation}%
\begin{equation}
B_{\Delta}(u,v)=\left\{
\begin{array}
[c]{ll}%
\int F(dx)~\sin(wv\Delta^{1-1/\beta}x), & \mbox{if }~\beta<1\\[3mm]%
\int F(dx)\left(  \sin(wv\Delta^{1-1/\beta}x)-wv\Delta^{1-1/\beta
}x1_{\{|x|\leq1\}}\right)  \qquad & \mbox{if }~\beta\geq1.
\end{array}
\right.  \label{eq:Ex5}%
\end{equation}
So we can inject these formulas directly into (\ref{eq:char1}).

As for the asymptotic variance in Theorem \ref{theorem:Gknown}, it is even
simpler. Indeed, we have here
\begin{equation}
\Psi_{k}(u,0)=e^{-w^{\beta}u^{\beta}/2}.\label{eq:char2}%
\end{equation}
Therefore $I(k)=-\Psi_{k}^{\prime}(1,0)=\beta~w^{\beta}~e^{-w^{\beta}/2}/2>0$
and $J(k)=\frac{1}{2}(\Psi_{k}(2,0)+1)-\psi_{k}(1,0)^{2}=\frac{1}{2}\left(
1+e^{-(2w)^{\beta}/2}\right)  -e^{-w^{\beta}}$, and thus
\begin{equation}
\Sigma^{2}(k)=2~\frac{1+e^{-(2w)^{\beta}/2}-2e^{-w^{\beta}}}{\beta
^{2}~w^{2\beta}~e^{-w^{\beta}}}~.\label{eq:Ex2}%
\end{equation}
When $\beta<2$, it turns out that the minimal variance is achieved for some
value $w=w_{\beta}\in(0,\infty)$, whereas $\Sigma^{2}(k)$ tends to $\infty$
when $w$ goes either to $0$ or to $\infty$. In contrast, when $\beta=2$ the
variance $\Sigma^{2}(k)$ goes to $1/2$ as $w\rightarrow0$: recall once more
that $1/2$ is the efficient\ variance in that case.

For the semiparametric situation, things are even simpler. The estimating
function of (\ref{FE6'}) becomes
\begin{equation}
U_{n,G,\beta,k}(u)=V_{n}(w)-\Psi_{k}\left(  \frac{u}{S_{n}},0\right)
,\label{eq:char10}%
\end{equation}
provided in (\ref{eq:charf}) we sum over $i\in\{q_{n}+m_{n}+1,\ldots,n\}$.
Moreover $u\mapsto\Psi_{k}(u,0)$ is invertible, so the estimator
$\widehat{\sigma}_{n}(k)$ takes the simple explicit form
\begin{equation}
\widehat{\sigma}_{n}(k)=S_{n}~\frac{2^{1/\beta}}{w}~\left(  -\log\left(
\frac{1}{p_{n}}\sum_{i=q_{n}+m_{n}+1}^{n}\cos\left(  \frac{w\Delta
_{n}^{-1/\beta}(\chi_{i}^{n}-B_{n})}{S_{n}}\right)  \right)  \right)
^{1/\beta}\label{eq:explicit}%
\end{equation}
if the argument of the logarithm is positive (otherwise, put for example
$\widehat{\sigma}_{n}(k)=1$).

\section{Example:\ Power and truncated power functions.\label{sec:Example2}}

Another natural choice for the function $k$ is a power function, that is
$k(x)=|x|^{r}$, for some $r>0$ when $\beta=2$ and $r\in(0,\beta/2)$ otherwise
(when $\beta=2$ this is -- in principle -- optimal for $r=2$). In general, the
function $\Psi_{G,\Delta_{n},\beta,k_{n}}$ is not explicit but can be
numerically approximated via Monte--Carlo procedures for example. We can also
compute the limiting variance: with the notation $m_{r}=E(|W_{1}|^{r})$ we get
$I(k)=-rm_{r}$ and $J(k)=m_{2r}-m_{r}^{2}$, hence
\begin{equation}
\Sigma^{2}(k)=\frac{m_{2r}-m_{r}^{2}}{r^{2}m_{r}^{2}}.\label{eq:Ex1}%
\end{equation}
When $\beta=2$ we have a closed expression for $m_{r}$ (see
(\ref{eq:PV_notrunc_nojump}) below), and not surprisingly $\Sigma^{2}(k)$
achieves its minimum, equal to $1/2$, at $r=2$: recall that $1/2$ is the
\textquotedblleft efficient\textquotedblright\ variance in that case. When
$\beta<2$ we have no explicit expression for these moments. However,
$\Sigma^{2}(k)$ goes to $\infty$ when $r$ increases to $\beta/2$, and we
conjecture that $\Sigma^{2}(k)$ is monotone increasing in $r$ (this property
holds at least when $\beta=1$); so one should take $r$ as small as possible,
although $r=0$ is of course excluded. \smallskip

In the semiparametric setting, the previous choice is not admissible, since
$k$ has to be bounded. So we must \textquotedblleft truncate\textquotedblright%
\ the argument, by using the following function $k=k_{\gamma}$:
\begin{equation}
k_{\gamma}(x)=|x|^{r}1_{\{|x|\leq\gamma\}}\label{eq:trunc}%
\end{equation}
for some constant $\gamma$. The function $\Psi_{k_{\gamma}}(u,0)=u^{r}%
E(|W_{1}|^{r}1_{\{|W_{1}|\leq\gamma/u\}})$ is invertible from a neighborhood
$I$ of $u=1$ onto some interval $I^{\prime}$, and we write $\Psi_{h_{\gamma}%
}^{-1}(v)$ for the inverse function at $v\in I^{\prime}$. Then if $B_{n}$ and
$S_{n}$ are the preliminary estimators, and if
\begin{equation}
V_{n}(\gamma)=\frac{1}{p_{n}\Delta_{n}^{r/\beta}}\sum_{i=m_{n}+1}^{n}|\chi
_{i}^{n}-B_{n}|^{r}1_{\{|\chi_{i}^{n}|\leq\gamma\Delta^{1/\beta}%
\}},\label{eq:vari}%
\end{equation}
the estimator $\widehat{\sigma}_{n}(k_{\gamma})$ is defined by
\begin{equation}
\widehat{\sigma}_{n}(k_{\gamma})=S_{n}~\Psi_{k_{\gamma}}^{-1}\left(
\frac{V_{n}(\gamma S_{n})}{S_{n}^{r}}\right)  \label{eq:vari2}%
\end{equation}
if the argument of $\Psi_{k_{\gamma}}^{-1}$ above is in $I^{\prime}$, and
$\widehat{\sigma}_{n}(k_{\gamma})=1$ (for example) otherwise. This is
almost\ as explicit as (\ref{eq:explicit}) is. Since $k_{\gamma}$ is even we
again have $J(k_{\gamma})=0$, whereas
\begin{equation}
\Sigma^{2}(k_{\gamma})=\frac{M_{\gamma,2r}-M_{\gamma,r}^{2}}{\left(
rM_{\gamma,r}-2h_{\beta}(\gamma)\gamma^{r+1}\right)  ^{2}},\quad
\mbox{where }~M_{\gamma,s}=\mathbf{E}(|W_{1}|^{r}1_{\{|W_{1}|\leq\gamma
\}}).\label{eq:var3}%
\end{equation}
We can then try to minimize this variance, by appropriately choosing the two
constants $\gamma>0$ and $r>0$.

One could also use $k_{\gamma_{n}}$, the $r$th power truncated at some level
$\gamma_{n}>0$ depending on $n$: our general results do not apply, but similar
results, with possibly other rates, should obviously apply. In fact, in the
next section we work out completely this kind of truncated power functions in
a particular case, to check that it is best (for the rate of convergence at
least) to take a constant level $\gamma_{n}=\gamma$, as it is implicitly
proposed in the method previously developed.

\section{Example:\ Brownian motion plus Gaussian compound Poisson
process.\label{sec:Example3}}

In this section, we present a fully worked out example, where $W$ is Brownian
motion and $Y\ $is a compound Poisson process with Gaussian jumps, say
$N(0,\eta)$, and intensity of jumps given by some $\lambda>0$. \cite{yacjfe04}
and \cite{mancini01} studied the estimation of the parameters of this model,
using a variety of methods.

As usual, we are interested in estimating the parameter $\sigma$ given the
increments $\chi_{i}^{n}$ of $X_{t}=\sigma W_{t}+Y_{t}$ (see (\ref{eq:Q1})).
We consider a number of estimating equations for this model, based on the
power or truncated power variations
\begin{equation}
V_{n}(c,\kappa)=\frac{1}{p_{n}\Delta_{n}^{r/2}}\sum_{i=m_{n}+1}^{n}|\chi
_{i}^{n}|^{r}1_{\{|\chi_{i}^{n}|\leq\tau(\Delta_{n})\}},\label{eq:trunc2}%
\end{equation}
for $r\in(0,2]$. Here $\tau(\Delta)$ is the truncation rate, taken to be of
the form $\tau(\Delta)=c\Delta^{1/2+\kappa}$ with $c$ a constant and
$\kappa\in(-1/2,\infty)$.

Note that $V_{n}$ above is exactly $V_{n}(\gamma)$ of (\ref{eq:vari}) with
$\gamma=\Delta_{n}^{-1/2}\tau(\Delta_{n})$ (here $Y$ is symmetrical, so
$B_{n}=0$). The associated estimator is then given by
\begin{equation}
\widehat{\sigma}_{n}=S_{n}~H_{\Delta_{n}}^{-1}\left(  \frac{V_{n}%
(cS_{n},\kappa)}{S_{n}^{r}}\right)  \label{eq:vari3}%
\end{equation}
where $H_{\Delta}^{-1}$ is the local inverse around $1$ of the function
$H_{\Delta}(u)=E(|uW_{\Delta}|^{r}1_{\{|uW_{\Delta}|\leq\tau(\Delta)\}})$.

When $c=\infty$ we get the (non truncated) $r$th power variation. If
$c<\infty$ and $\kappa=0$ this corresponds to taking $k=k_{c}$, as given by
(\ref{eq:trunc}): we essentially eliminate from the sum above the increments
in which $Y$ jumps. When $\kappa>0$ we eliminate more increments, and fewer
when $\kappa<0$.

The expected values of the powers without truncation are given by
\begin{equation}
\mathbf{E}(|X_{\Delta}|^{r})=\sum_{j=0}^{+\infty}\frac{2^{r/2}}{\sqrt{\pi}%
j!}~\Gamma\left(  \frac{1+r}{2}\right)  ~e^{-\lambda\Delta}~(\lambda
\Delta)^{j}\left(  \sigma^{2}\Delta+j\eta\right)  ^{r/2},
\label{eq:PV_notrunc}%
\end{equation}%
\begin{equation}
\mathbf{E}(|\sigma W_{\Delta}|^{r})=\frac{2^{r/2}}{\sqrt{\pi}}~\Gamma\left(
\frac{1+r}{2}\right)  \sigma^{r}\Delta^{r/2} \label{eq:PV_notrunc_nojump}%
\end{equation}
With truncation at rate $\tau(\Delta),$ we get
\begin{align}
&  \mathbf{E}\left(  |X_{\Delta}|^{r}1_{\{|X_{\Delta}|\leq\tau(\Delta
)\}}\right) \label{eq:PV_trunc}\\
&  \quad=e^{-\lambda\Delta}\sum_{j=0}^{+\infty}\frac{2^{r/2}}{\sqrt{\pi}%
j!}\left(  \Gamma\left(  \frac{1+r}{2}\right)  -\Gamma\left(  \frac{1+r}%
{2},\frac{\tau(\Delta)^{2}}{2\left(  \sigma^{2}\Delta+j\eta\right)  }\right)
\right)  (\lambda\Delta)^{j}\left(  \sigma^{2}\Delta+j\eta\right)
^{r/2}\nonumber
\end{align}
where $\Gamma(a,\cdot)$ denotes the incomplete Gamma function of order $a$,
and
\begin{equation}
\mathbf{E}\left(  |\sigma W_{\Delta}|^{r}1_{\{|\sigma W_{\Delta}|\leq
\tau(\Delta)\}}\right)  =\frac{2^{r/2}}{\sqrt{\pi}}\left(  \Gamma\left(
\frac{1+r}{2}\right)  -\Gamma\left(  \frac{1+r}{2},\frac{\tau(\Delta)^{2}%
}{2\sigma^{2}\Delta}\right)  \right)  \sigma^{r}\Delta^{r/2}.
\label{eq:PV_trunc_nojump}%
\end{equation}
When $r=2,$ we have $\Gamma\left(  3/2\right)  =\sqrt{\pi}/2$ and
$\Gamma\left(  \frac{3}{2},x\right)  =e^{-x}\sqrt{x}+\sqrt{\pi}~\Phi(\sqrt
{2x})$ where $\Phi$ denotes the cdf of the $N(0,1)$ law. Similarly simpler
expressions are also obtained in the case where $r=1,$ since $\Gamma\left(
1\right)  =1$ and $\Gamma\left(  1,x\right)  =e^{-x}.$

As described above, in the semiparametric case where the distribution of $Y$
is not known to the statistician, we propose to use an approximate centering
based on computing these expectations assuming that $X=\sigma W$ only (i.e.,
as if there were no jumps) and we will study the behavior of this estimator
when $Y$ is in fact a compound Poisson process. The effect of the
misspecification error is to bias the resulting estimator of $\sigma.$ But, at
the leading order in $\Delta,$ the expected values of the moments functions
computed without jumps coincide with those computed under the correct
specification. Indeed, for $X$ from (\ref{eq:X=sW+aY}), we have
\begin{align*}
\mathbf{E}\left(  |X_{\Delta}|^{r}\right)   &  =\mathbf{E}\left(  |\sigma
W_{\Delta}|^{r}\right)  +o(\Delta^{r/2})\\
\mathbf{E}\left(  |X_{\Delta}|1_{\{|X_{\Delta}|\leq\tau(\Delta)\}}\right)   &
=\mathbf{E}\left(  |\sigma W_{\Delta}|1_{\{|\sigma W_{\Delta}|\leq\tau
(\Delta)\}}\right)  +o(\Delta^{r/2}),
\end{align*}
with the second result following from
\begin{equation}
\Gamma(a,x)~=~\left\{
\begin{array}
[c]{ll}%
\Gamma(a)+x^{a}\left(  -\frac{1}{a}+\frac{x}{1+a}+O(x^{2})\right)  \qquad &
\text{ \ \ near }~0\\[2.5mm]%
e^{-x}x^{-1+a}\left(  1+\frac{a-1}{x}+O(x^{-2})\right)  \qquad & \text{
\ \ near }+\infty.
\end{array}
\right.  \label{eq:Gamma}%
\end{equation}

As a result, the bias of the estimator of $\sigma$ based on approximate
centering will vanish asymptotically in $\Delta$ and we will have a result of
the form
\[
\sqrt{n\Delta_{n}^{v_{1}}}~\left(  \widehat{\sigma}_{n}-\bar{\sigma}%
_{n}\right)  \rightarrow N\left(  0,v_{0}\right)
\]
where
\[
\bar{\sigma}_{n}=\sigma+b_{0}\Delta_{n}^{b_{1}}+o(\Delta_{n}^{b_{1}})
\]
with $b_{1}>0$. (If $b_{1}=0$ for some choice of $(r,\kappa,c)$ then the
parameter $\sigma$ is not identified by an estimating function based on that
combination.) Also, $v_{1}=0$ corresponds to a rate of convergence of the
estimator of $n^{1/2},$ and any value $v_{1}>0$ corresponds to a slower than
$n^{1/2}$ rate of convergence.

We also note that when $b_{1}>0$ the rate of convergence and asymptotic
variance of the semiparametric estimator of $\sigma$ are identical at the
leading order in $\Delta_{n}$ to the expressions one would obtain in the fully
parametric, correctly specified, case where centering of the estimating
equation is done with either (\ref{eq:PV_notrunc}) or (\ref{eq:PV_trunc}) as
appropriate, instead of the approximate centering using
(\ref{eq:PV_notrunc_nojump}) or (\ref{eq:PV_trunc_nojump}). Centering using
the latter is of course the only feasible estimator in the semiparametric case
where the distribution of $Y$ is unknown.

In what follows, we use the explicitness of this model to fully characterize
the asymptotic distribution of the semiparametric estimator of $\sigma$, i.e.,
$(b_{0},b_{1},v_{0},v_{1})$ as functions of $(r,\kappa,c)$ and the parameters
of the model $(\sigma,\lambda,\eta).$

\subsection{Power variations without truncation.}

In that situation, we have for the asymptotic variance:

\begin{itemize}
\item When $0<r<1,$ we have $v_{1}=0$ and $v_{0}=\frac{1}{r^{2}}\left(
\sqrt{\pi}\frac{\Gamma\left(  \frac{1}{2}+r\right)  }{\Gamma\left(  \frac
{1+r}{2}\right)  ^{2}}-1\right)  $.

\item When $r=1,$ we have $v_{1}=0$ and $v_{0}=\frac{1}{2}~\left(  \left(
\pi-2\right)  \sigma^{2}+\pi\lambda\eta\right)  $.

\item When $1<r<2,$ we have $v_{1}=r-1$ and $v_{0}=\frac{\sqrt{\pi}%
\sigma^{2-2r}\lambda\eta^{r}}{r^{2}}~\frac{\Gamma\left(  \frac{1}{2}+r\right)
}{\Gamma\left(  \frac{1+r}{2}\right)  ^{2}}$.
\end{itemize}

As for the bias, when $0<r<2$ we have $b_{1}=1-r/2$ and $b_{0}=\frac
{\sigma^{1-r}\lambda\eta^{r/2}}{r}$.

\begin{remark}
The estimator based on power variations converges (not taking the bias into
consideration) at rate $n^{1/2}$ only when $r\leq1.$ When $r>1$ the mixture of
jumps and volatility slows down the rate of convergence ($v_{1}>0$). When
$r=2,$ the parameter $\sigma$ is simply not identified, as is obvious from the
fact that $E(X_{\Delta}^{2})=$ $(\sigma^{2}+\lambda\eta)\Delta.$ This is also
apparent here from the fact that $b_{1}\downarrow0$ as $r\uparrow2,$ so the
bias no longer vanishes asymptotically. And the bias even worsens the rate, of course.
\end{remark}

\begin{remark}
When $r<1,$ the asymptotic variance $v_{0}$ is identical to the expression
obtained without jumps, as was the case when the log-likelihood score was used
as an estimating equation. When $r=1,$ the rate of convergence remains
$n^{1/2},$ but $v_{0}$ is larger in the presence of jumps.
\end{remark}

\subsection{Power variations with $\Delta^{1/2}$ truncation.}

If we truncate the increments according to $\tau(\Delta)=c\Delta^{1/2},$ then
$v_{1}=0$ for all values of $r\in(0,2]$ and
\[
v_{0}=\frac{2^{r}\sigma^{4+2r}\left(  \sqrt{\pi}\left(  \Gamma\left(  \frac
{1}{2}+r\right)  -\Gamma\left(  \frac{1}{2}+r,\frac{c^{2}}{2\sigma^{2}%
}\right)  \right)  -\left(  \Gamma\left(  \frac{1+r}{2}\right)  -\Gamma\left(
\frac{1+r}{2},\frac{c^{2}}{2\sigma^{2}}\right)  \right)  ^{2}\right)
}{\left(  \sqrt{2}c^{1+r}\exp\left(  -\frac{c^{2}}{2\sigma^{2}}\right)
-2^{r/2}r\sigma^{1+r}\left(  \Gamma\left(  \frac{1+r}{2}\right)
-\Gamma\left(  \frac{1+r}{2},\frac{c^{2}}{2\sigma^{2}}\right)  \right)
\right)  ^{2}}%
\]

As for the bias, we have $b_{1}=1$ and
\[
b_{0}=\frac{\sigma\lambda\left(  \Gamma\left(  \frac{1+r}{2}\right)
-\Gamma\left(  \frac{1+r}{2},\frac{c^{2}}{2\sigma^{2}}\right)  \right)
}{\left(  \Gamma\left(  \frac{1+r}{2}\right)  -\Gamma\left(  \frac{1+r}%
{2},\frac{c^{2}}{2\sigma^{2}}\right)  \right)  -2\left(  \Gamma\left(
\frac{3+r}{2}\right)  -\Gamma\left(  \frac{3+r}{2},\frac{c^{2}}{2\sigma^{2}%
}\right)  \right)  }.
\]

\begin{remark}
Truncating at rate $\Delta^{1/2}$ restores the convergence rate $n^{1/2}$ for
all values of $r,$ (again, regardless of the bias) and permits identification
when $r=2.$ When $0<r<1$ (where the rate $n^{1/2}$ was already achieved
without truncation), not truncating can lead to either a smaller or larger
value of $v_{0}$ than truncating at rate $n^{1/2},$ depending upon the values
of $(\sigma^{2},c).$
\end{remark}

\begin{remark}
The asymptotic variance $v_{0}$ is identical to its expression when no jumps
are present, as it should be in view of our general results (as said before,
this type of truncation leads to the estimators studied in our general
results). In all cases, the bias is smaller than when no truncation is applied.
\end{remark}

\subsection{Power variations with slower than $\Delta^{1/2}$ truncation.}

If we now keep too many increments by truncating according to $\tau
(\Delta)=c\Delta^{1/2+\kappa},$ with $-1/2<\kappa<0,$ then we have for
$r\in(0,2]:$

\begin{itemize}
\item When $-3/(2+4r)<\kappa<0,$ we have $v_{1}=0$ and
\[
v_{0}=\frac{\sigma^{2}}{r^{2}}\left(  \sqrt{\pi}\frac{\Gamma\left(  \frac
{1}{2}+r\right)  }{\Gamma\left(  \frac{1+r}{2}\right)  ^{2}}-1\right)
\]

\item When $\kappa=-3/(2+4r),$ we have $v_{1}=0$ and
\[
v_{0}=\frac{2^{1/2-r}c^{1+2r}\sqrt{\pi}\lambda\sigma^{2-2r}}{r^{2}\left(
1+2r\right)  \eta^{1/2}\Gamma\left(  \frac{1+r}{2}\right)  ^{2}}+\frac
{\sigma^{2}}{r^{2}}\left(  \sqrt{\pi}\frac{\Gamma\left(  \frac{1}{2}+r\right)
}{\Gamma\left(  \frac{1+r}{2}\right)  ^{2}}-1\right)
\]

\item When $-1/2<\kappa<-3/(2+4r),$ we have $v_{1}=-\kappa-2r\kappa-3/2>0$
and
\[
v_{0}=\frac{2^{1/2-r}c^{1+2r}\sqrt{\pi}\lambda\sigma^{2-2r}}{r^{2}\left(
1+2r\right)  \eta^{1/2}\Gamma\left(  \frac{1+r}{2}\right)  ^{2}}.
\]

\end{itemize}

As for the bias, we have:

\begin{itemize}
\item When $-1/(2+2r)<\kappa<0,$ we have $b_{1}=1$ and $b_{0}=-\frac
{\lambda\sigma}{r}$

\item When $\kappa=-1/(2+2r),$ we have $b_{1}=1$ and $b_{0}=\frac
{\lambda\sigma}{\left(  1+r\right)  }\left(  \frac{2^{1/2-r/2}c^{1+r}}%
{r\sqrt{\eta}\sigma^{r}\,\Gamma\left(  \frac{1+r}{2}\right)  }-1-\frac{1}%
{r}\right)  $

\item When $-1/2<\kappa<-1/(2+2r),$ we have $b_{1}=3/2+\kappa+r\kappa>0$ and
$b_{0}=\frac{2^{1/2-r/2}c^{1+r}\lambda\sigma^{1-r}}{r\,\left(  1+r\right)
\sqrt{\eta}\Gamma\left(  \frac{1+r}{2}\right)  }$.
\end{itemize}

\begin{remark}
When $0<r<1,$ we are automatically in the situation where $\kappa>-3/(2+4r)$,
and hence keeping more than $O(\Delta_{n}^{1/2})$ increments results in the
convergence rate $n^{1/2}$ and the same asymptotic variance $v_{0}$ as when
keeping all increments (i.e., not truncating at all). When $1<r<2,$ however,
it is possible to restore the convergence rate $n^{1/2}$ (compared to not
truncating)\ by keeping more than $O(\Delta_{n}^{1/2})$ increments, but still
\textquotedblleft not too many\textquotedblright\ of them ($-3/(2+4r)\leq
\kappa<0$) beyond that; but even keeping a larger fraction of the increments
($-1/2<\kappa<-3/(2+4r)$) results in an improvement over keeping all
increments since $3/2-\kappa-2r\kappa<r-1$ so that the rate of convergence of
$\widehat{\sigma}_{n},$ although slower than $n^{1/2},$ is nonetheless faster
than $n^{1/2}\Delta_{n}^{(r-1)/2}.$
\end{remark}

\begin{remark}
The expressions for $\kappa<0$ do not converge to those with $O(\Delta
_{n}^{1/2})$ truncation as $\kappa\uparrow0$ because of the essential
singularity of the incomplete $\Gamma$ function near infinity, given in
(\ref{eq:Gamma}):\ when $\tau(\Delta)=c\Delta^{1/2+\kappa}$ then
$\Gamma((1+r)/2,\cdot)$ is evaluated at $\tau(\Delta)^{2}/(2\sigma^{2}\Delta)$
$=$ $c^{2}\Delta^{2\kappa}/(2\sigma^{2})$ and for fixed $\kappa<0,$ terms
proportional to $\exp(-c^{2}\Delta^{2\kappa}/(2\sigma^{2}))$ are negligible in
the Taylor series in $\Delta$ of $v_{0}$ and $b_{0}$. This is not the case
when $\kappa=0$ however.
\end{remark}

\begin{remark}
As for the bias, keeping \textquotedblleft too many\textquotedblright\ but not
all increments ($-1/2<\kappa<-1/(2+2r)$) leads to a smaller bias than keeping
all increments, since $3/2+\kappa+r\kappa$ $>$ $1-r/2,$ but to a larger bias
than keeping just the right amount since $3/2+\kappa+r\kappa<1.$
\end{remark}

\subsection{Power variations with faster than $\Delta^{1/2}$ truncation.}

Finally, if we keep too few increments by truncating according to
$c\Delta^{1/2+\kappa},$ with $\kappa>0,$ then $v_{1}=\kappa$ for all values of
$r\in(0,2]$ and
\[
v_{0}=\frac{\sqrt{2\pi}\left(  1+r\right)  ^{2}\sigma^{3}}{2c\left(
1+2\,r\right)  }%
\]
As for the bias, we have $b_{1}=1$ and $b_{0}=\sigma\lambda$.

\begin{remark}
Truncating at a rate faster than $\Delta^{1/2}$ deteriorates the convergence
rate of the estimator from $n^{1/2}$ to $n^{1/2}\Delta_{n}^{\kappa/2}$:\ while
we successfully eliminate the impact of jumps on the estimator, we are at the
same time reducing the effective sample size utilized to compute the
estimator, which increases its asymptotic variance.
\end{remark}

\begin{remark}
The expressions for $v_{0}$ and $b_{0}$ for $\kappa>0$ also do not converge to
those with $O(\Delta_{n}^{1/2})$ truncation as $\kappa\downarrow0$ because
once again we cannot interchange the order of the limits $\Delta
_{n}\rightarrow0$ and $\kappa\rightarrow0$.
\end{remark}

\subsection{Comparison with the general case.}

Let us compare, in the semiparametric case, the specific results just obtained
with the general results obtained in Theorems \ref{theorem:Gunknown1} and
\ref{theorem:Gunknown}. In the present situation we have $G\in\overline
{\mathcal{G}}_{0}^{\prime}$. So these general results assert that if
\begin{equation}
n\Delta_{n}^{2}\rightarrow0,\label{eq:aseff0}%
\end{equation}
then the estimators $\widehat{\sigma}_{n}$ converge at a rate $\sqrt{n}$, and
the limit of the normalized error is Gaussian without bias; when
(\ref{eq:aseff0}) fails but $\Delta_{n}\rightarrow0$ yet, then the sequence
$((\sqrt{n}\bigwedge\Delta_{n}^{-1})(\widehat{\sigma}_{n}-\sigma)$ is tight.

The estimators (\ref{eq:vari3}) converge at rate $\sqrt{n}$ when $v_{1}=0$ and
$n\Delta_{n}^{2b_{1}}$ is bounded (then there is a bias) or $n\Delta
_{n}^{2b_{1}}\rightarrow0$ (there is no bias). Otherwise, the sequence
$(\sqrt{n\Delta_{n}^{v_{1}}}\bigwedge\Delta_{n}^{-b_{1}})(\widehat{\sigma}%
_{n}-\sigma)$ is tight. Then:

\begin{itemize}
\item Power variation without truncation: we have a rate $\sqrt{n}$ only when
$r\in(0,1]$ and $n\Delta_{n}^{2-r}$ is bounded. Otherwise the rate is always
worse than in our general results: this was expected, of course.

\item Power variation with $\Delta^{1/2}$ truncation: \ If $n\Delta_{n}%
^{2}\rightarrow0$ we have rate $\sqrt{n}$ with asymptotically unbiased error.
If $n\Delta_{n}^{2}\rightarrow a\in(0,\infty)$ we have rate $\sqrt{n}$ with
asymptotically biased error. If $n\Delta_{n}^{2}\rightarrow\infty$, then
$\Delta_{n}^{-1}(\widehat{\sigma}_{n}-\sigma)$ converges in probability to the
constant $b_{0}$: this is a bit better than what we get by applying the
general results recalled above. This holds irrespectively of $r\in(0,2]$ (and
also for $r>2$ here, as a matter of fact), but of course the asymptotic
variance depends on $r$, and also on $c$.

\item Power variation with slower than $\Delta^{1/2}$ truncation: \ The rate
is $\sqrt{n}$ if $-1/(2+2r)\leq\kappa<0$ and $n\Delta_{n}^{2}$ is bounded, or
if $-3(2+4r)\leq\kappa<-1/(2+2r)$ and $n\Delta_{n}^{3+2\kappa+2r\kappa}$ is
bounded. This is worse than the previous case.

\item Power variation with faster than $\Delta^{1/2}$ truncation: \ The rate
is at most $\sqrt{n\Delta_{n}^{\kappa}}$, and always worst than in the
$\Delta^{1/2}$ truncation case.
\end{itemize}

\section{Example: Sum of two stable processes.\label{sec:Example4}}

In this last section we consider the case where $Y$ is also a symmetric stable
process, with index $\alpha\in(0,\beta)$. Then $G\in\overline{\mathcal{G}%
}_{\alpha}^{\prime}$.

\subsection{The empirical characteristic function.}

First, we can consider estimators based on the empirical characteristic
function, that is we consider $k(x)=\cos(wx)$ for some $w>0$. We have the
parametric estimate $\widehat{\sigma}_{n}=\widehat{\sigma}_{n}(G,\phi,k)$ of
Theorem \ref{theorem:Gknown} (here $k$ is bounded, so $\phi$ is indeed
irrelevant). The sequence $\sqrt{n}~(\widehat{\sigma}_{n}-\sigma)$ converges
in law to $N(0,\sigma^{2}\Sigma^{2}(k))$, where $\Sigma^{2}(k)$ is given by
(\ref{eq:Ex2}). On the other hand we have the semiparametric estimators
$\widehat{\sigma}_{n}(k)$, which by Theorem \ref{theorem:Gunknown} behaves as
such: under
\begin{equation}
n\Delta_{n}^{\frac{2(\beta-\alpha)}{\beta}}\rightarrow0, \label{eq:aseff0'}%
\end{equation}
$\sqrt{n}~(\widehat{\sigma}_{n}(k)-\sigma)$ converges in law to $N(0,\sigma
^{2}\Sigma^{2}(k))$. And in general the sequence $(\sqrt{n}\bigwedge\Delta
_{n}^{-\frac{\beta-\alpha}{\beta}})(\widehat{\sigma}_{n}-\sigma)$ is tight.

In fact, since we are in Case 2 the preliminary estimator $S_{n}=S_{n}(G)$ is
the same in both cases, and $\widehat{\sigma}_{n}$ and $\widehat{\sigma}%
_{n}(k)$ are the solution of $U_{n}(u)=0$ and $U_{n}^{\prime}(u)=0$
respectively, which are closest to $S_{n}$, and the difference between these
two estimating functions is
\[
U_{n}(u)-U_{n}^{\prime}(u)=\widehat{U}_{n}(u):=\Psi_{G,\Delta_{n},\beta
,k}\left(  \frac{u}{S_{n}},\frac{1}{S_{n}},0\right)  -\Psi_{k}\left(  \frac
{u}{S_{n}},0\right)
\]
(recall (\ref{eq:char1}) and (\ref{eq:char10})). If we use the explicit forms
(\ref{eq:Ex3}) and (\ref{eq:char2}), we get
\[
\widehat{U}_{n}(u)=e^{-w^{\beta}u^{\beta}/2S_{n}^{b}}~\left(  e^{w^{\alpha
}\Delta_{n}^{\frac{\beta-\alpha}{\beta}}/2S_{n}\alpha}-1\right)  ,
\]
which is equivalent to $\frac{w^{\alpha}}{2\sigma^{\alpha}}~\Delta_{n}%
^{\frac{\beta-\alpha}{\beta}}~e^{-w\beta/2}$ as $n\rightarrow\infty$ and
$u\rightarrow\sigma$ (recall that $S_{n}\rightarrow\sigma$ in probability).
Since $\Psi_{k}^{\prime}(1,0)=-\beta e^{-w^{b}/2}\neq0$, we deduce that the
difference $\widehat{\sigma}_{n}(k)-\widehat{\sigma}_{n}$ is equivalent (in
probability) to $-(w^{a}/2\beta\sigma^{\alpha})\Delta_{n}^{\frac{\beta-\alpha
}{\beta}}$. Therefore, in addition to the fact that $\sqrt{n}~(\widehat
{\sigma}_{n}(k)-\sigma)$ converges in law to $N(0,\sigma^{2}\Sigma^{2}(k))$
under (\ref{eq:aseff0'}), we get

\begin{itemize}
\item If $n\Delta_{n}^{\frac{\beta-\alpha}{\beta}}\rightarrow a^{2}%
\in(0,\infty)$, then $\sqrt{n}~(\widehat{\sigma}_{n}(k)-\sigma)$ converges in
law to $N(-aw^{a}/2\beta\sigma^{\alpha},\sigma^{2}\Sigma^{2}(k))$,

\item If $n\Delta_{n}^{\frac{\beta-\alpha}{\beta}}\rightarrow\infty$, then
$\Delta_{n}^{-\frac{\beta-\alpha}{\beta}}~(\widehat{\sigma}_{n}(k)-\sigma)$
converges in probability to the constant $-w^{a}/2\beta\sigma^{\alpha}$.
\end{itemize}

We conclude that the results of Theorem \ref{theorem:Gunknown} are sharp, for
the particular estimation functions $k(x)=\cos(wx)$ at least.

\subsection{Truncated power functions.}

We can do a similar analysis for the estimators (\ref{eq:vari2}), based on the
truncated power variation $V_{n}(\gamma)$ of (\ref{eq:vari}) with $B_{n}=0$
(because $Y$ is symmetrical here). That is, we consider the truncated power
variations at the level $\Delta_{n}^{1/\beta}$. Namely when $n\Delta
_{n}^{2\frac{\beta-\alpha}{\beta}}\rightarrow\infty$, one can show that, at
least when $\gamma$ is small enough (but it is probably true for all $\gamma>0
$), then the sequence $\Delta_{n}^{-\frac{\beta-\alpha}{\beta}}(\widehat
{\sigma}_{n}-\sigma)$ is tight and its limiting distributions include some
Dirac masses at non vanishing constants. So here again the results of Theorem
\ref{theorem:Gunknown} are sharp. But of course, as already said before, this
does not completely rule out the existence of estimators constructed in a
different way and behaving better.

\section{Conclusions.}

We exhibited a class of estimators for the volatility parameter $\sigma$ in a
model where the driving process $W_{t}$ is perturbed by another process
$Y_{t}.$ These estimators can be designed in such a way that they are immune
to the presence of the perturbation $Y_{t}:$ they are asymptotically
efficient, in the strong\ sense that they behave asymptotically like the
efficient estimators for the model $X_{t}=\sigma W_{t}$ with no perturbing
term $Y_{t}$.\bigskip

\bigskip
%

\renewcommand\baselinestretch{1.1}
\normalsize

\bigskip%

\Line{\AOSaddress{{Department of Economics}\\
{Princeton University and NBER}\\
{Princeton, NJ 08544-1021}\\
{E-mail: yacine@princeton.edu}}\hfill\AOSaddress
{{Laboratoire de Probabilit\'es (UMR 7599)}\\
{Universit\'e P. et M. Curie (Paris-6)}\\
{75252 Paris C\'edex 05 }\\
{E-mail: jj@ccr.jussieu.fr}}}%

\end{document}